# PREHOMOGENEOUS VECTOR SPACES EISENSTEIN SERIES AND INVARIANT THEORY

Akihiko Yukie[1]

Oklahoma State University

ABSTRACT. In this paper, we give an introduction to the rationality of the equivariant Morse stratification, and state the author's results on zeta functions of prehomogeneous vector spaces.

**Content**
   Introduction
   §1 Review of invariant theory
   §2 Examples of Morse stratifications
   §3 An invariant measure on $\mathrm{GL}(n)$
   §4 The zeta function
   §5 The smoothed Eisenstein series
   §6 Adjusting terms
   §7 The general program
   §8 Examples of principal part formulas
   References

**Introduction**

The aim of this paper is to state the author's results on the global theory of zeta functions associated with prehomogeneous vector spaces and give an introduction to equivariant Morse stratifications and their rationality.

Let $\mathbb{Z}, \mathbb{Q}, \mathbb{R}, \mathbb{C}$ be the set of integers, rational numbers, real numbers, and complex numbers respectively. Let $k$ be a number field, and $o_k$ its integer ring. Let $\mathbb{A}$ and $\mathbb{A}^\times$ be the sets of adeles and ideles of $k$ respectively. We use the notation $|x|$ for the adelic absolute value of $x \in \mathbb{A}$. Also let $|x|_\infty$ be the product of absolute values of $x$ at all the infinite places.

Let $G$ be a connected reductive group, $V$ a representation of $G$, and $\chi_V$ an indivisible non-trivial character of $G$ all defined over $k$. In this paper, we are mainly interested in the case when $G$ is a product of $\mathrm{GL}(n)$'s. The triple $(G, V, \chi_V)$ is called a prehomogeneous vector space if the following two conditions are satisfied.

  (1) There exists a Zariski open orbit.
  (2) There exists a polynomial $0 \neq \Delta \in k[V]$ such that $\Delta(gx) = \chi'(g)\Delta(x)$ where $\chi'$ is a rational character, proportional to $\chi_V$.

[1]Partially supported by NSF Grant DMS-9101091

A polynomial $\Delta$ which satisfies the property (2) is called a relative invariant polynomial. We define $V_k^{\text{ss}} = \{x \in V_k \mid \Delta(x) \neq 0\}$, and call it the set of semi-stable points. For practical purposes, we assume that there exists a one dimensional split torus $T_0$ contained in the center of $G$ such that if $\alpha \in T_0$, it acts on $V$ by multiplication by $\alpha$ and $\chi_V(\alpha) = \alpha^e$ for some $e > 0$.

Our purpose in this paper is to outline the author's program to compute Laurent expansions at poles of zeta functions for prehomogeneous vector spaces when the group is a product of $\operatorname{GL}(n)$'s.

The notion of prehomogeneous vector spaces and their zeta functions was introduced by Sato and Shintani in [17]. However, some special cases were studied earlier by many mathematicians, most notably Siegel [20] and Davenport [4]. Neither Siegel nor Davenport used the zeta function explicitly to obtain the density theorems for equivalence classes of quadratic forms [20] and binary cubic forms [4]. Davenport's result was later used by Davenport and Heilbronn [5], [6] to prove the density theorem for cubic extensions of $\mathbb{Q}$, which states

$$\sum_{\substack{[k:\mathbb{Q}]=3 \\ |\Delta_k|_\infty \leq x}} 1 \sim \frac{x}{\zeta(3)},$$

where $k$ runs through all the cubic fields and $\Delta_k$ is the discriminant of $k$.

D. J. Wright and the author proved in [23] that generic $G_k$-orbits of several "parabolic type" prehomogeneous vector spaces correspond almost bijectively with number field extensions of certain degrees. Therefore, by studying the global and the local zeta functions for those cases, we can expect certain density theorems where number fields are counted possibly with some weights.

Let $L \subset V_k^{\text{ss}}$ be a $G_k$-invariant subset, and $\Phi$ a Schwartz–Bruhat function on $V_\mathbb{A}$. Let $\omega$ be a character of $G_\mathbb{A}/G_k$. Consider an integral of the following form

$$(0.1) \qquad Z_L(\Phi, \omega, s) = \int_{G_\mathbb{A}/G_k} |\chi_V(g)|^s \omega(g) \sum_{x \in L} \Phi(gx) dg,$$

where $dg$ is a Haar measure on $G_\mathbb{A}$. If (0.1) is well defined for $L = V_k^{\text{ss}}$ and $\operatorname{Re}(s)$ sufficiently large, we say that $(G, V, \chi_V)$ is of complete type. Otherwise, we say that $(G, V, \chi_V)$ is of incomplete type. If $(G, V, \chi_V)$ is a prehomogeneous vector space of complete type, the zeta function is, by definition, the integral $Z_{V_k^{\text{ss}}}(\Phi, \omega, s)$.

Suppose that $(G, V, \chi_V)$ is a prehomogeneous vector space of complete type. Then the meromorphic continuation and the functional equation of the zeta function were proved by reducing the problem to the local theory at infinite places (see [17], [15]). Therefore, in most cases, the Tauberian theorem implies the existence of a limit of the form

$$\sum_{\substack{y \in G_{o_k} \backslash V_{o_k}^{\text{ss}} \\ |P(y)|_\infty \leq x}} \mu_y \sim cx^a (\log x)^b,$$

where $\mu_y$ is a number determined by the stabilizer of $y$, and $c, a, b$ are constants. The improved version of the Tauberian theorem in [17] gives us the error term estimate also. However, the existence of a similar limit for $G_k$-orbits is unknown for most cases and does not automatically follow from the the meromorphic continuation and the functional equation of the zeta function, or the knowledge of the location and the orders of the poles. The process to handle $G_k$-orbits of prehomogeneous vector spaces from the viewpoint of zeta function theory is called the "filtering process" and was formulated by Datskovsky and Wright [3]. This process has to be modified for the general case, because the zeta function can have higher order poles in general, and it is likely that the possible generalization of this process requires answers to the following two problems.

Problem A Determine the Laurent expansions at all the poles of the zeta function.

Problem B Compute the local zeta function explicitly for special test functions.

If we use the filtering process, the only pole we can use is the rightmost pole. However, it is likely that all the poles have to be studied in order to say anything about the rightmost pole. The author's approach to Problem A is a generalization of Shintani's approach in [18], and is based on the framework of equivariant Morse theory. We give an introduction to invariant theory and related materials for the sake of number theorists in §1,2.

If $f(s)$ is a meromorphic function having the Laurent expansion $\sum a_i(s-s_0)^i$ at $s=s_0$, we call $a_i$'s generalized special values of $f$ at $s=s_0$. Roughly speaking, we want to express the coefficients of the Laurent expansions at the poles of the zeta function by the generalized special values of the zeta functions which correspond to "unstable strata" in equivariant Morse theory.

Let $G_{\mathbb{A}}^1 = \{g \in G_{\mathbb{A}} \mid |\chi_V(g)| = 1\}$. It is easy to prove (see [17]) that if $(G, V, \chi_V)$ is of complete type and the function $\sum_{x \in V_k} \Phi(gx)$ is integrable on $G_{\mathbb{A}}^1/G_k$, all the poles of the zeta function are simple and the residues can be computed. However, this assumption hardly applies to any important cases. Even though Siegel [20] and Davenport [4] essentially computed the residue of the zeta function at the rightmost pole for their cases, Shintani's work on the space of binary cubic forms is the first non-trivial work where Problem A was answered from the viewpoint of zeta function theory.

The primary difficulty of Problem A is based on the difficulty of handling divergent integrals. So far, our main analytic tool to handle this difficulty is the smoothed Eisenstein series and the iterated use of the Poisson summation formula. We discuss the smoothed Eisenstein series in §5 and describe our program in §7.

Among the eight cases in [23], four types $(C_2, D_4, D_5, E_6)$ are of incomplete type. Also even if we start with a prehomogeneous vector space of complete type, we have to face prehomogeneous vector spaces of incomplete type in the induction step anyway. In that case, we have to use generalized special values of the "adjusted zeta functions" (see [26]). For this reason, understanding of prehomogeneous vector spaces of incomplete type is essential to the entire theory. Therefore, from the viewpoint of Problem A, we can not restrict ourselves to prehomogeneous vector spaces of complete type.

There are not so many prehomogeneous vector spaces of incomplete type for which the Laurent expansions at the poles of the zeta function is known. The space of binary quadratic forms is an example of such a case. The difficulty of handling prehomogeneous vector spaces of incomplete type is that some semi-stable points have extra stabilizing elements. As is shown in [19], [25], [26] the zeta function itself does not satisfy a functional equation. Instead, it is natural to consider an adjusted version of the zeta function. We discuss this issue in §6.

In §8, we list some examples of prehomogeneous vector spaces for which Problem A was answered.

This manuscript was prepared mainly during the academic year 1990–1991 when the author was staying at SFB 170 Göttingen. The author would like to thank SFB 170 for supporting this project.

**§1 Review of invariant theory**

In this section, we review invariant theory and equivariant Morse theory. We do not restrict ourselves to prehomogeneous representations in this section.

Let $k$ be a perfect field. Let $G$ be a split connected reductive group and $V$ a representation of $G$ both defined over $k$. This $G$ does *not* correspond to the group $G$ in the introduction. Instead, it corresponds to $G' = \text{Ker}(\chi_V)$. Note that since $G$ is connected and $\chi_V$ is indivisible, $G'$ is connected.

For any group $G$ over $k$, let $X^*(G), X_*(G)$ be the groups of rational characters, rational one parameter subgroups (which we abbreviate to "1PS" from now on) of $G$ respectively. We choose a maximal split torus $T \subset G$. We define $\mathfrak{t} = X_*(T) \otimes \mathbb{R}$, $\mathfrak{t}^* = X^*(T) \otimes \mathbb{R}$. Let $\mathfrak{t}_\mathbb{Q} = X_*(T) \otimes \mathbb{Q}$, $\mathfrak{t}_\mathbb{C} = X_*(T) \otimes \mathbb{C}$ etc. Elements in $\mathfrak{t}_\mathbb{Q}, \mathfrak{t}_\mathbb{Q}^*$ are called rational elements. We fix a Weyl group invariant inner product $(\ ,\ )$ on $\mathfrak{t}^*$. We assume that if $z, z' \in \mathfrak{t}_\mathbb{Q}^*$, then $(z, z') \in \mathbb{Q}$. Let $\|\ \|$ be the metric defined by this inner product. We choose a split Borel subgroup $T \subset B$ and a Weyl chamber $\mathfrak{t}_-^* \subset \mathfrak{t}^*$ so that the weights of the unipotent radical of $B$ belong to $\mathfrak{t}_-^*$ by the conjugation $b \to tbt^{-1}$. We can identify $\mathfrak{t}, \mathfrak{t}^*$ by this inner product, and therefore, we consider $\|\ \|, \mathfrak{t}_-$ for $\mathfrak{t}$ also.

We recall the definition of stability over $\bar{k}$. Let $\pi : V \setminus \{0\} \to \mathbb{P}(V)$ be the natural projection map. Let $\bar{k}[V]^{G_{\bar{k}}}$ be the ring of polynomials invariant under the action of $G_{\bar{k}}$. Suppose that $P \in \bar{k}[V]^{G_{\bar{k}}}$ is a homogeneous polynomial. We define $\mathbb{P}(V)_P = \{\pi(x) \mid P(x) \neq 0\}$.

**Definition (1.1)** *Let $y \in \mathbb{P}(V)_{\bar{k}}$. We define*
*(1) $y$ is semi-stable if there exists a homogeneous polynomial $P \in \bar{k}[V]^{G_{\bar{k}}}$ such that $y \in \mathbb{P}(V)_P$,*
*(2) $y$ is properly stable if there exists a homogeneous polynomial $P \in \bar{k}[V]^{G_{\bar{k}}}$ such that $y \in \mathbb{P}(V)_P$, all the orbits in $\mathbb{P}(V)_P$ are closed, and the stabilizer of $y$ in $G_{\bar{k}}$ is finite,*
*(3) $y$ is unstable if it is not semi-stable.*

We use the notation $\mathbb{P}(V)^{\text{ss}}_{\bar{k}}$ (resp. $\mathbb{P}(V)^s_{(0)\bar{k}}$) for the set of semi-stable (resp. properly stable) points. These are $\text{Gal}(\bar{k}/k)$-invariant opensets of $\mathbb{P}(V)_{\bar{k}}$. However, we will consider the set $\mathbb{P}(V)^s_{(0)\bar{k}}$ in number theoretic situations only when $\mathbb{P}(V)^{\text{ss}}_{\bar{k}} = \mathbb{P}(V)^s_{(0)\bar{k}}$.

We will recall three equivalent definitions of the equivariant Morse stratification for the rest of this section. The first definition is invariant theoretic and is closely related to the Hilbert–Mumford criterion of stability. We need this definition, because we are interested in the rationality question in number theoretic situations. This approach is due to Kempf and Ness. The second definition is combinatorial and was given by both Kirwan and Ness. We need this definition for practical determination of the equivariant Morse stratification. The third definition is Morse theoretic and is based on the convexity property of the moment map. The moment map and the Morse theoretic approach was used by Kirwan to prove the smoothness and the perfectness of the stratification when $k = \mathbb{C}$.

(1) The invariant theoretic definition.

Let $k$ be an algebraically closed field. We are only interested in the plain stability, so our situation corresponds to the case $S = \{0\}$ in [8]. Let $X$ be a variety and $f : \mathbb{G}_m \to X$ a morphism. We imbed $\mathbb{G}_m$ to the one dimensional affine space $\mathrm{Aff}^1$ in the usual manner. We say that $\lim_{\alpha \to 0} f(\alpha) = y$ if $f$ extends to a morphism from $\mathrm{Aff}^1$ to $X$ and $f(0) = y$. If $\lambda$ is a 1PS of $G$, we define $P(\lambda) = \{p \in G \mid \lim_{\alpha \to 0} \lambda(\alpha) p \lambda(\alpha)^{-1}$ exists $\}$. The group $P(\lambda)$ is a parabolic subgroup of $G$.

Let $x \in V \setminus \{0\}$. We define $|V, x|$ (resp. $|V, x|_{\{0\}}$) to be the set of 1PS's $\lambda$ such that $\lim_{\alpha \to 0} \lambda(\alpha) x$ exists (resp. $\lim_{\alpha \to 0} \lambda(\alpha) x = 0$). Let $\lambda$ be a non-trivial 1PS of $G$ (not necessarily in $|V, x|$). We consider when $\lambda$ belongs to $|V, x|_{\{0\}}$. Let $x = \sum_{i=1}^{N} x_i$ be the eigen decomposition with respect to $\lambda$, i.e. $\lambda(\alpha) x = \sum_{i=1}^{N} \alpha^{p_i} x_i$, $x_i \neq 0$ for all $i$, and $p_i \neq p_j$ if $i \neq j$.

**Definition (1.2)**

(1) $$\mu(x, \lambda) = \min_{i=1}^{N} p_i,$$

(2) $$\nu(x, \lambda) = \frac{\mu(x, \lambda)}{\|\lambda\|}.$$

If $\lambda_1, \lambda_2$ are two non-trivial proportional 1PS's, then $\nu(x, \lambda_1) = \nu(x, \lambda_2)$.

The following theorem is Theorem (3.4) [8] ((1),(2) were proved by Mumford in Proposition 2.17 [12]).

**Theorem (1.3)**
(1) *The function $\nu(x, \lambda)$ has a maximum value (which we denote by $B_x$) on the set $|V, x|$ if it is not empty.*
(2) *The condition $\overline{Gx} \ni 0$ is equivalent to the condition that $|V, x|_{\{0\}}$ is non-empty.*

*Suppose that $|V, x|_{\{0\}}$ is non-empty, and let $\Lambda_x$ be the set of indivisible 1PS $\lambda$'s such that $\nu(\lambda, x) = B_x$.*
(3) *The set $\Lambda_x$ is non-empty, and there exits a parabolic subgroup $P_x$ of $G$ such that $P_x = P(\lambda)$ for all $\lambda \in \Lambda_x$.*
(4) *The set $\Lambda_x$ is a principal homogeneous space under the action of the unipotent radical of $P_x$.*
(5) *Any maximal torus of $P_x$ contains a unique element of $\Lambda_x$.*

We followed [8] for the definition of $\mu(x, \lambda)$. However, this is $-\mu^L(x, \lambda)$ in [12]. Roughly speaking, if $x \in V$ corresponds to a geometric object, elements of $\Lambda_x$ are 1PS's which give the fastest degeneration of $x$. Now we are ready to state the stratification defined by Ness [13].

Suppose that $x \in V \setminus \{0\}$, and $\pi(x)$ is unstable. This means that $\Lambda_x$ is non-empty. There is a unique 1PS $\lambda \in \mathfrak{t}_-$ such that $\lambda$ is conjugate to any element of $\Lambda_x$. Let $\lambda' \in \Lambda_x$. Let $\beta_x \in \mathfrak{t}_-$ be the element which is proportional to $\lambda$ and $\|\beta_x\| = \nu(x, \lambda')$. This $\beta_x$ is uniquely determined by $x$, and is a rational element.

For $\beta \in \mathfrak{t}_- \cap \mathfrak{t}_\mathbb{Q}$, we define

(1.4) $$S_\beta = \{x \in V \setminus \{0\} \mid \beta_x = \beta\}.$$

Let $S_0 = \pi^{-1}(\mathbb{P}(V)^{\mathrm{ss}})$. Apparently,

$$V \setminus \{0\} = \coprod_\beta S_\beta,$$

where $\beta$ runs through all the elements of $\mathfrak{t}_- \cap \mathfrak{t}_\mathbb{Q}$. It turns out that we only have to consider finitely many such $\beta$'s. We discuss this issue in the second definition. In [13], Ness used the notion of the lengths of vectors (see [9]) to formulate the stratification. The fact that the formulation in [9] is equivalent to the above formulation is explained in Lemma 3.1 [10].

Hilbert-Mumford criterion of stability says that a point is semi-stable or properly stable if and only if it is stable or properly stable with respect to the action of all the non-trivial 1PS's. The above theorem implies that if a point is unstable, there is a canonical choice of a 1PS which, in a way, gives the worst degeneration of $x$. Since we are interested in number theoretic applications, one natural question is if we can choose such a 1PS which is split over the ground field. This question was answered affirmatively by Kempf [8].

The following theorem was proved by Kempf (Theorem (4.2) [8]).

**Theorem (1.5)** *Suppose that $x \in V_k$, and $|V, x|_{\{0\}}$ is non-empty over the closure $\bar{k}$. Then all the elements in $\Lambda_x$ are split, $P_x$ is rationally conjugate to a standard parabolic subgroup, and $\Lambda_x$ is a principal homogeneous space under the action of the $k$-points of the unipotent radical of $P_x$.*

The above theorem was proved without the assumption that the group is split. This theorem implies that if $x \in V_k$ and $|V, x|_{\{0\}}$ is non-empty over the closure $\bar{k}$, $x$ is conjugate to $y \in V_k$ by an element of $G_k$ and $\Lambda_y$ contains a split 1PS of $T$. We discuss the inductive structure of $S_\beta$ and its rationality in the next definition.

(2) The combinatorial definition

Let $k, G, V$ etc. be as before. We can diagonalize the action of $T$ because $k$ is a perfect field. Let $(x_0, \cdots, x_N)$ be a coordinate system of $V$ whose coordinate vectors are eigenvectors of $T$. Let $\gamma_i \in \mathfrak{t}^*$ be the weight determined by the $i$-th coordinate. The following definition is due to Kirwan.

**Definition (1.6)** *A point $\beta \in \mathfrak{t}_\mathbb{Q}^*$ is called a minimal combination of weights if $\beta$ is the closest point to the origin of the convex hull of a finite subset of $\{\gamma_1, \cdots, \gamma_N\}$.*

We denote the set of minimal combination of weights which lie in $\mathfrak{t}_-^*$ by $\mathfrak{B}$. In this definition, $\mathfrak{B}$ is the index set of the stratification. Since $\mathfrak{B}$ is a finite set, this means that the stratification is a finite stratification. It turns out that $S_0$ is the set of semi-stable points. We describe the stratum $S_\beta$ for $\beta \in \mathfrak{B} \setminus \{0\}$. The reader should note that $S_\beta$ can be the empty set.

We change the ordering of $\{\gamma_1, \cdots, \gamma_N\} = \{\gamma_1', \cdots, \gamma_N'\}$ so that

$$(\gamma_1', \beta) \leq \cdots \leq (\gamma_N', \beta).$$

We assume that

$$((\gamma'_1, \beta), \cdots, (\gamma'_N, \beta)) = \frac{1}{m_0}(\overbrace{m_1, \cdots, m_1}^{n_1}, \overbrace{m_2, \cdots, m_2}^{n_2}, \cdots, \overbrace{m_p, \cdots, m_p}^{n_p})$$

where $m_0 > 0, m_1 < \cdots < m_p$ are coprime integers. Since $\beta \in \mathfrak{B}$, there exists some $1 \leq s \leq p$ such that $\frac{m_s}{m_0} = \|\beta\|^2$. Let $e'_i$ be the coordinate vector which corresponds to $\gamma'_i$. Let $V_b$ be the subspace spanned by

$$\{e'_i \mid (\gamma'_i, \beta) = \frac{m_b}{m_0}\}.$$

We define

$$Z_\beta = \oplus_{b=s} V_b, W_\beta = \oplus_{b>s} V_b, Y_\beta = Z_\beta \oplus W_\beta,$$
$$\bar{Z}_\beta = \{\pi(x) \mid x \in Z_\beta \setminus \{0\}\}, \bar{Y}_\beta = \{\pi(x,y) \mid x \in Z_\beta \setminus \{0\}, y \in W_\beta\},$$
$$M_\beta = \{g \in G \mid Ad(g)\beta = \beta\}.$$

The group $M_\beta$ is $\text{Stab}_\beta$ in [10].

Let $p_\beta : \bar{Y}_\beta \to \bar{Z}_\beta$ be the projection map. Let $\lambda_\beta$ be the indivisible rational character of $M_\beta$ whose restriction to $T$ is a positive multiple of $\beta$. Let $\nu_\beta$ be the indivisible split 1PS of $M_\beta$ which is a positive multiple to $\beta$. We define $G_\beta = \{g \in M_\beta \mid \lambda_\beta(g) = 1\}$. Note that $G_\beta$ is connected. The group $G_\beta$ acts on $\bar{Z}_\beta$ linearly. Let $\bar{Z}_\beta^{ss}$ be the set of semi-stable points with respect to this action. Let $P_\beta$ be the standard parabolic subgroup of $G$ whose Levi component is $M_\beta$ and fixes the set $Y_\beta$. Let $U_\beta$ be the unipotent radical of $P_\beta$. Since $\lambda_\beta$ is split, all these groups are split groups over $k$.

Let $\bar{Y}_{\beta\bar{k}}^{ss} = p_\beta^{-1}(\bar{Z}_{\beta\bar{k}}^{ss})$. We define $S_{\beta\bar{k}} = \pi^{-1}(\bar{S}_{\beta\bar{k}})$ etc. The following inductive structure of the strata was proved by Kirwan and Ness [10], [13].

**Proposition (1.7)** $\bar{S}_{\beta\bar{k}} \cong G_{\bar{k}} \times_{P_{\beta\bar{k}}} \bar{Y}_{\beta\bar{k}}^{ss}$.

The fact that this $S_\beta$ and $S_\beta$ in (1) coincide is proved in [10, pp. 150–156]. More precisely, suppose that $x \in S_{\beta\bar{k}}$, and $\lambda \in \Lambda_x$. Choose $g \in G_{\bar{k}}$ so that $g\lambda g^{-1}$ is a 1PS of $T_{\bar{k}}$ proportional to $\beta$. Then $gx$ must belong to $Y_{\bar{k}}^{ss}$.

(3) The convexity property of the moment map

We assume that $k = \mathbb{C}$. Let $K$ be a maximal compact subgroup of $G$, and $T_K$ a maximal torus of $K$. We assume that the action of $T_K$ on $\mathbb{P}(V)$ is diagonalized. Let $\mathfrak{g}, \mathfrak{k}, \mathfrak{t}$, be the Lie algebras of $G, K, T_K$ respectively. Let $\mathfrak{t}_-$ be a Weyl chamber. We choose a $K$ invariant Hermitian product $(\ ,\ )_V$ on $V$. This defines a homomorphism of Lie algebras $\mathfrak{k} \mapsto \mathfrak{u}(N)$ (the $N$ dimensional unitary Lie algebra). We consider the Fubini-Study metric on $\mathbb{P}(V)$ defined by this metric. Let $\omega$ be its Kahler form. If $a \in \mathfrak{k}$, it defines a vector field on $\mathbb{P}(V)$.

**Definition (1.8)** *A map $\mu : \mathbb{P}(V) \mapsto \mathfrak{k}^*$ is called a moment map if it satisfies the following two conditions*
*(1) $\mu$ is $K$-equivariant with respect to the co-adjoint action of $K$ on $\mathfrak{k}^*$,*
*(2) $d\mu(x)(\xi)(a) = \omega_x(\xi, a_x)$ for all $x \in \mathbb{P}(V)$, $\xi \in T_x\mathbb{P}(V)$, and $a \in \mathfrak{k}$.*

If a moment map exists and $K$ is semi-simple, it is uniquely determined by the above conditions.

For $\mathbb{P}(V)$, one can construct a moment map in the following manner. Suppose that $x \in V \setminus \{0\}$, $\pi(x) = y \in \mathbb{P}(V)$, and $a \in \mathfrak{u}(N)$. We define a map $\bar{\mu} : \mathbb{P}(V) \to \mathfrak{u}(N)^*$ by

$$\bar{\mu}(y)(a) = \frac{{}^t\bar{x}ax}{2\pi\sqrt{-1}(x,x)}.$$

Let $\mu$ be the composition of $\bar{\mu}$ and the natural map $\mathfrak{u}(N)^* \mapsto \mathfrak{k}^*$. This $\mu$ is a moment map.

Let $\mu_T$ be the composition of $\mu$, and the projection $\mathfrak{k}^* \to \mathfrak{t}^*$. The third definition of the stratification is based on the following convexity property of the moment map (see [1], [7] or, [12] Appendix by Mumford).

**Theorem (1.9)** *Let $X \subset \mathbb{P}(V)$ be a connected $G$-invariant closed subvariety. Then $\mu_T(X) \cap \mathfrak{t}_-^*$ is a rational convex polytope.*

For $y = \pi(x) \in \mathbb{P}(V)$, we consider the closest point of $\mu(\overline{Gy}) \cap \mathfrak{t}_-^*$ to the origin, which we call $\beta_x$. Let $\beta \in \mathfrak{t}_{\mathbb{Q}}^* \cap \mathfrak{t}_-^*$. By the identification $\mathfrak{t} \cong \mathfrak{t}^*$, we can consider $\beta$ as an element of $\mathfrak{t}_{\mathbb{Q}} \cap \mathfrak{t}_-$. We define $S_\beta = \{x \mid \beta_x = \beta\}$. Clearly,

$$V \setminus \{0\} = \coprod_\beta S_\beta.$$

This definition coincides with previous definitions. Kirwan proved that this stratification is perfect with respect to $H_G(\ ) \otimes \mathbb{Q}$ (the equivariant cohomology group) using this formulation. The isomorphism $H_G^*(\mathbb{P}(V)) \otimes \mathbb{Q} \cong H^*(BG) \otimes H^*(\mathbb{P}(V)) \otimes \mathbb{Q}$ is proved in 5.8 [10] (also see [2]). Also

$$H_G^*(\bar{S}_\beta) \otimes \mathbb{Q} \cong H_{P_\beta}^*(\bar{Y}_\beta^{\mathrm{ss}}) \otimes \mathbb{Q} \cong H_{M_\beta}^*(\bar{Z}_\beta^{\mathrm{ss}}) \otimes \mathbb{Q} \cong H^*(B\mathbb{C}^\times) \otimes H_{G_\beta}^*(\bar{Z}_\beta^{\mathrm{ss}}) \otimes \mathbb{Q}.$$

Therefore, this gives an effective way of computing $H_G^*(\mathbb{P}(V)^{\mathrm{ss}}) \otimes \mathbb{Q}$. The basic reason why the stratification is perfect is that by $\lambda_\beta$, $\mathbb{C}^\times$ acts non-trivially on $Z_\beta$ by scalar multiplications, and it induces a multiplication by an element of $H^*(B\mathbb{C}^\times) \otimes H_{G_\beta}^*(\bar{Z}_\beta^{\mathrm{ss}}) \otimes \mathbb{Q}$ whose first component is non-zero via the Gysin map (see Theorem 5.4 [10] for the details).

This is the Morse stratification with respect to the function $\|\mu\|^2$. However, critical points are not isolated, so one has to give the direction of the Morse flow, which is given by $\lambda_\beta$ for points in $Y_\beta^{\mathrm{ss}}$.

Now we consider the situation in the introduction. Let $(G, V, \chi_V)$ be a prehomogeneous vector space over a number field $k$. Let $G' = \mathrm{Ker}(\chi_V)$ as before. The group $G'$ is connected as we mentioned at the beginning of this section. Let

$$V_k \setminus \{0\} = \coprod_\beta S_{\beta k}$$

be the equivariant Morse stratification with respect to the action of $G'$.

**Lemma (1.10)** $G_k = G'_k T_k$.

*Proof.* We only have to show that $\chi_V$ is indivisible as a character of $T$ also. Since all the characters of $G$ are rational, we assume that $k$ is algebraically closed. Let $C$ be the connected center, and $\bar G = [G, G]$. Then $G = C\bar G$, $\bar G$ is semi-simple and $\bar T = T \cap \bar G$ is a maximal torus of $\bar G$. Suppose that $\chi_V = \chi^p$ on $T$. Since $\bar T$ is connected, this implies that $\chi$ is trivial on $\bar T$ also. Since $T/\bar T \cong C/C \cap \bar G \cong G/\bar G$, we can extend $\chi$ to a character of $G$. Regular elements are dense in $G$, so $\chi_V = \chi^p$. Therefore, $p = \pm 1$.

Q.E.D.

Let $M_\beta = \{g \in G \mid Ad_g(\beta) = \beta\}$, and $M'_\beta = M_\beta \cap G'$ etc. Kempf's theorem (1.5) and (1.7) imply the following proposition.

**Proposition (1.11)** $S_{\beta k} \cong G_k \times_{P_{\beta k}} Y^{ss}_{\beta k}$.

*Proof.* Let $x \in S_{\beta k}$. Then any $\lambda \in \Lambda_x$ is split. Therefore, we can choose $g$ in the remark after (1.7) from $G_k$. So $gx \in Y^{ss}_{\bar k} \cap V_k = Y^{ss}_k$. Therefore, $G_k \times Y^{ss}_k \to S_{\beta k}$ is surjective. Suppose that $g_1, g_2 \in G'_k$, $t_1, t_2 \in T_k$, $y_1, y_2 \in Y^{ss}_k$, and $g_1 t_1 y_1 = g_2 t_2 y_2$. Then $g_1 \pi(t_1 y_1) = g_2 \pi(t_2 y_2)$. By (1.8), there exists $p \in P'_{\beta \bar k}$ such that $g_1 = pg_2$, $\pi(t_1 y_1) = p^{-1} \pi(t_2 y_2)$. From the first equation, $p \in P'_{\beta k}$. By assumption, $g_2 p t_1 y_1 = g_2 t_2 y_2$. Therefore, $y_1 = t_1^{-1} p^{-1} t_2 y_2$. Let $\widetilde p = t_2^{-1} p t_1$. This $\widetilde p$ is an element of $P_{\beta k}$, and satisfies the condition $g_1 t_1 = g_2 t_2 \widetilde p$, $y_1 = \widetilde p^{-1} y_2$.

Q.E.D.

Equivariant Morse theory in algebraic situation like in this section was established around 1983. However, Kempf's paper appeared in Annals in 1978. Therefore, he proved the rationality of the Morse stratification before equivariant Morse theory was established.

Of course, if $(G, V, \chi_V)$ is a prehomogeneous vector space over an algebraically closed field, the equivariant Morse stratification is the decomposition into $G$-orbits. But how do we determine the orbit decompositions of prehomogeneous vector spaces systematically? It is easy to observe from the second definition of the stratification that it is possible to write a computer program to determine all the possible $\mathfrak{B} \ni \beta$'s. Moreover, the knowledge of $\mathfrak{B}$ enables us to determine the inductive structure of each stratum also.

## §2 Examples of Morse stratifications

We consider equivariant Morse stratifications for the following examples
(1) $G = \mathrm{GL}(1) \times \mathrm{GL}(2)$, $V = \mathrm{Sym}^3 k^2$,
(2) $G = \mathrm{GL}(2) \times \mathrm{GL}(1)^2$, $V = \mathrm{Sym}^2 k^2 \oplus k^2$,
(3) $G = \mathrm{GL}(3) \times \mathrm{GL}(2)$, $V = \mathrm{Sym}^2 k^3 \otimes k^2$.

The zeta function for the case (1) was studied by Shintani. The case (2) is an example of a reducible representation and was studied by F. Sato [14]. The case (3) is the $F_4$ case (or the quartic case) in [23].

For later purposes, we fix some notations. Let

(2.1) $\quad \mathfrak{t}_{A_{n-1}} = \{q = (q_1, \cdots, q_n) \in \mathbb{R}^n \mid q_1 + \cdots + q_n = 0\}.$

Let $\mathfrak{t}_{A_{n-1}\mathbb{C}}$ be the complexification of $\mathfrak{t}_{A_{n-1}}$. For

$$q = (q_1, \cdots, q_n), q' = (q'_1, \cdots, q'_n) \in \mathfrak{t}_{A_{n-1}},$$

we define $(q, q')_{A_{n-1}} = \sum_i q_i q'_i$. For $z = (z_1, \cdots, z_n) \in \mathfrak{t}_{A_{n-1}\mathbb{C}}$ and a permutation $\tau$, we define $\tau z = (z_{\tau(1)}, \cdots, z_{\tau(n)})$. This inner product is invariant under permutations. By this inner product, we identify $\mathfrak{t}_{A_{n-1}}$ with its dual space. Let $\|\ \|_{A_{n-1}}$ be the metric defined by this inner product.

(1) The case $G = \mathrm{GL}(1) \times \mathrm{GL}(2)$, $V = \mathrm{Sym}^3 k^2$

The group $\mathrm{GL}(2)$ acts on $V$ naturally, and if $\alpha \in \mathrm{GL}(1)$, $\alpha$ acts by multiplication by $\alpha$. Let $\widetilde{T} \subset G$ be the kernel of the homomorphism $G \to \mathrm{GL}(V)$. We consider the character $\chi_V(t, g) = (\det g)^2 t^{-3}$ for $(t, g) \in G$. We can consider $\chi_V$ as a character of $G/\widetilde{T}$, and it is indivisible as a character of $G$. It is known that $(G/\widetilde{T}, V, \chi_V)$ is a prehomogeneous vector space. Let $G'$ be the kernel of $\chi_V$. Then $G'$ is isomorphic to $\mathrm{SL}(2)$. Therefore, $\mathfrak{t}, \mathfrak{t}^*$ can be identified with $\mathfrak{t}_{A_1}$. We choose an identification $\mathbb{R} \cong \mathfrak{t}^*$ by the map $q \to (-q, q)$. We choose $\mathfrak{t}^*_- = \{q \in \mathbb{R} \mid q \geq 0\}$ as the Weyl chamber. The weighs of coordinates of $V$ are $-3, -1, 1, 3$. Therefore, there are two possibilities for $\beta$, i.e. $\beta_1 = 1, \beta_2 = 3$. The stratum $S_{\beta_1}$ is the set of forms with one double factor. The stratum $S_{\beta_2}$ is the set of forms with one triple factor (see [18], [22] for the details).

(2) The case $G = \mathrm{GL}(2) \times \mathrm{GL}(1)^2$, $V = \mathrm{Sym}^2 k^2 \oplus k^2$

The group $\mathrm{GL}(2)$ acts on each irreducible component of $V$ naturally. Let $g = (g_1, t_1, t_2)$ and $v = (v_1, v_2)$, where

$$g_1 \in \mathrm{GL}(2),\ t_1, t_2 \in \mathrm{GL}(1),\ v_1 \in \mathrm{Sym}^2 k^2, v_2 \in k^2.$$

We define $gv = (t_1 g_1 v_1, t_2 g_1 v_2)$. Let $\widetilde{T} \subset G$ be the kernel of the homomorphism $G \to \mathrm{GL}(V)$. Let $\chi_V(g) = (\det g_1)^a t_1^{b_1} t_2^{b_2}$ be a character of $G$. The character $\chi_V$ is trivial on $\widetilde{T}$ if and only if $2a - 2b_1 - b_2 = 0$. Suppose that this condition is satisfied and $(a, b_1, b_2)$ is relatively prime. Then $\chi_V$ is indivisible as a character of $G$. We consider $\chi_V$ as a character of $G/\widetilde{T}$. The kernel of $\chi_V$ is isomorphic to $\mathrm{SL}(2) \times \mathrm{GL}(1)$ by $(g_1, \alpha) \to (g_1, \alpha^{b_2}, \alpha^{-b_1})$. We assume that $2b_1 > b_2$. We identify $\mathfrak{t}^*$ with $\mathbb{R}^2$ by the map $q = (q_1, q_2) \to ((-q_1, q_1), b_2 q_2, -b_1 q_2)$. We define $((q_1, q_2), (q'_1, q'_2)) = 2 q_1 q'_1 + \frac{1}{14} q_2 q'_2$. We choose $\{(q_1, q_2) \in \mathbb{R}^2 \mid q_1 \geq 0\}$ as the Weyl chamber. The weights of the coordinates look as the following picture.

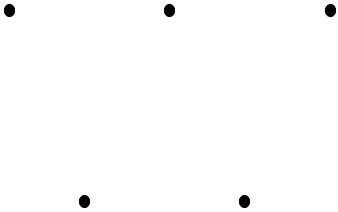

We call these points $\gamma_{1,11}, \gamma_{1,12}, \gamma_{1,22}, \gamma_{2,1}, \gamma_{2,2}$ from the left to the right and from the top to the bottom. Consider the convex hulls of the following sets
   (1) $C_1 : \{\gamma_{1,22}, \gamma_{2,1}\}$,
   (2) $C_2 : \{\gamma_{1,12}, \gamma_{2,2}\}$,
   (3) $C_3 : \{\gamma_{1,22}, \gamma_{2,2}\}$,
   (4) $C_4 : \{\gamma_{1,11}, \gamma_{1,22}\}$,

(5) $C_5 : \{\gamma_{1,22}\}$,
(6) $C_6 : \{\gamma_{2,2}\}$.

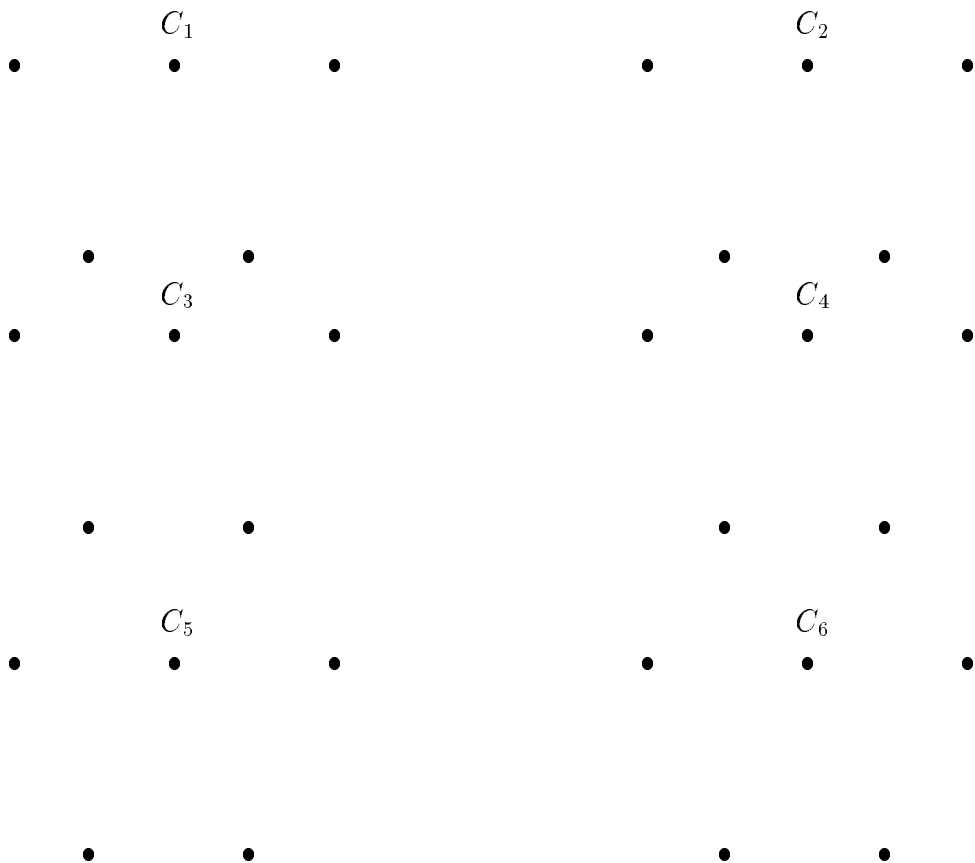

Let $\beta_i$ be the closet point to the origin from $C_i$ for $i = 1, \cdots, 6$. In this case, there are precisely 6 unstable strata $S_{\beta_1}, \cdots, S_{\beta_6}$ (The vector $(0, -b_1)$ gives the empty set).

(3) The case $G = \mathrm{GL}(3) \times \mathrm{GL}(2)$, $V = \mathrm{Sym}^2 k^3 \otimes k^2$

Let $\widetilde{V}$ be the space of quadratic forms in three variables $v = (v_1, v_2, v_3)$. We identify $\widetilde{V}$ with $k^6$ as follows

$$Q_x(v) = x_{11} v_1^2 + x_{12} v_1 v_2 + x_{13} v_1 v_3 + x_{22} v_2^2 + x_{23} v_2 v_3 + x_{33} v_3^2$$
$$\to x = (x_{11}, x_{12}, x_{13}, x_{22}, x_{23}, x_{33}).$$

The group $G_1 = \mathrm{GL}(3)$ acts on $\widetilde{V}$ in the following way

$$Q_{g_1 \cdot x}(v) = Q_x(v g_1),$$

for $g_1 \in G_1$.

We identify $V$ with $\widetilde{V} \oplus \widetilde{V}$. Any element of $V$ is of the form

$$Q = (Q_{x_1}, Q_{x_2}),$$

where
$$x_1 = (x_{1,11}, x_{1,12}, x_{1,13}, x_{1,22}, x_{1,23}, x_{1,33}),$$
$$x_2 = (x_{2,11}, x_{2,12}, x_{2,13}, x_{2,22}, x_{2,23}, x_{2,33}).$$

We use $x = (x_1, x_2)$ as the coordinate system of $V$.

If $g = (g_1, g_2) \in \mathrm{GL}(3) \times \mathrm{GL}(2)$ and $g_2 = \begin{pmatrix} a & b \\ c & d \end{pmatrix}$, we define

$$g \cdot Q = (aQ_{g_1 \cdot x_1} + bQ_{g_1 \cdot x_2}, cQ_{g_1 \cdot x_1} + dQ_{g_1 \cdot x_2}).$$

Let $\widetilde{T}$ be as before. Let $\chi_V(g_1, g_2) = (\det g_1)^4 (\det g_2)^3$ for $g_1 \in \mathrm{GL}(3), g_2 \in \mathrm{GL}(2)$. The character $\chi_V$ is tirvial on $\widetilde{T}$, and we consider $\chi_V$ as a character of $G/\widetilde{T}$. The character $\chi_V$ is indivisible as a character of $G$.

First, we fix a metric on $\mathfrak{t}^*$. We identify $\mathfrak{t}^*$ with $\mathfrak{t}^*_{A_2} \oplus \mathfrak{t}^*_{A_1}$. For $q = (q_1, q_2), q' = (q'_1 q'_2) \in \mathfrak{t}^*$, we define $(q, q') = (q_1, q'_1)_{A_2} + (q_2, q'_2)_{A_1}$. This inner product is Weyl group invariant.

We define

$$\mathfrak{t}^*_- = \{q = (q_{11}, q_{12}; q_{13}, q_{21}, q_{22}) \in \mathfrak{t}^* \mid q_{11} \leq q_{12} \leq q_{13},\ q_{21} \leq q_{22}\}.$$

We use $\mathfrak{t}^*_-$ as the Weyl chamber. The weights of $x_{1,ij}$, $x_{2,ij}$ are as follows

| | | | | |
|---|---|---|---|---|
| $x_{1,11}$ | $(\frac{4}{3}, -\frac{2}{3}, -\frac{2}{3}; \frac{1}{2}, -\frac{1}{2})$ | | $x_{2,11}$ | $(\frac{4}{3}, -\frac{2}{3}, -\frac{2}{3}; -\frac{1}{2}, \frac{1}{2})$ |
| $x_{1,12}$ | $(\frac{1}{3}, \frac{1}{3}, -\frac{2}{3}; \frac{1}{2}, -\frac{1}{2})$ | | $x_{2,12}$ | $(\frac{1}{3}, \frac{1}{3}, -\frac{2}{3}; -\frac{1}{2}, \frac{1}{2})$ |
| $x_{1,13}$ | $(\frac{1}{3}, -\frac{2}{3}, \frac{1}{3}; \frac{1}{2}, -\frac{1}{2})$ | | $x_{2,13}$ | $(\frac{1}{3}, -\frac{2}{3}, \frac{1}{3}; -\frac{1}{2}, \frac{1}{2})$ |
| $x_{1,22}$ | $(-\frac{2}{3}, \frac{4}{3}, -\frac{2}{3}; \frac{1}{2}, -\frac{1}{2})$ | | $x_{2,22}$ | $(-\frac{2}{3}, \frac{4}{3}, -\frac{2}{3}; -\frac{1}{2}, \frac{1}{2})$ |
| $x_{1,23}$ | $(-\frac{2}{3}, \frac{1}{3}, \frac{1}{3}; \frac{1}{2}, -\frac{1}{2})$ | | $x_{2,23}$ | $(-\frac{2}{3}, \frac{1}{3}, \frac{1}{3}; -\frac{1}{2}, \frac{1}{2})$ |
| $x_{1,33}$ | $(-\frac{2}{3}, -\frac{2}{3}, \frac{4}{3}; \frac{1}{2}, -\frac{1}{2})$ | | $x_{2,33}$ | $(-\frac{2}{3}, -\frac{2}{3}, \frac{4}{3}; -\frac{1}{2}, \frac{1}{2})$ |

With our metric, the weights of the coordinates look as the following picture.

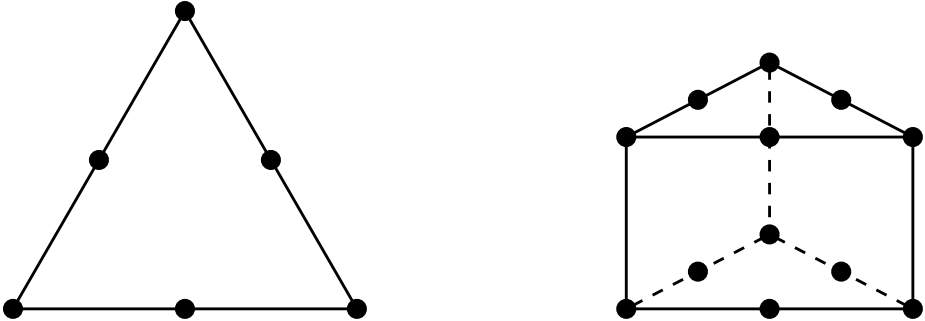

We proved in [23] that $q \in V_k^{\mathrm{ss}}$ if and only if the zero set of $Q_1 = Q_2 = 0$ in $\mathbb{P}^2_k$ consists of precisely four points. The unstable strata also have geometric meanings. Tables (2.2) (2.3) are the list of unstable strata in this case (see [26] for the details).

Table (2.2)

| Strata | Convex hull | Conics |
|---|---|---|
| $S_{\beta_1}$ | 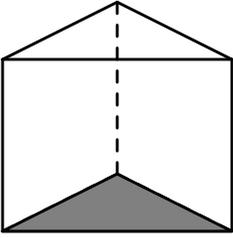 | 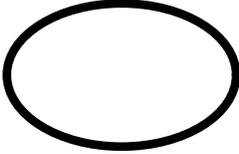<br>2 identical non-singular conics |
| $S_{\beta_2}$ | 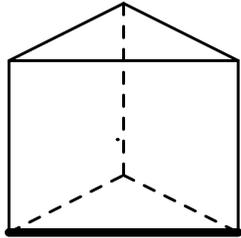 | 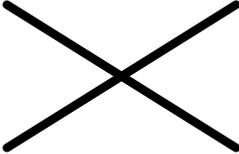<br>2 identical reducible conics |
| $S_{\beta_3}$ | 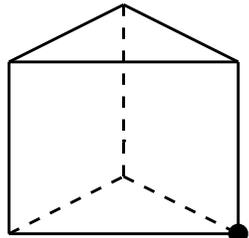 | 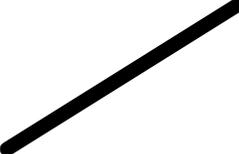<br>2 identical double lines |
| $S_{\beta_4}$ | 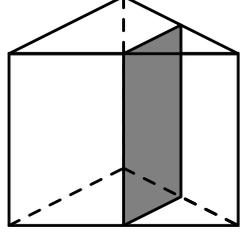 | 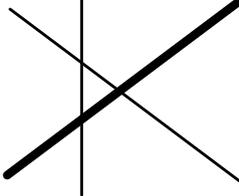<br>1 common component |

| Strata | Convex hull | Conics |
|---|---|---|
| $S_{\beta_5}$ | 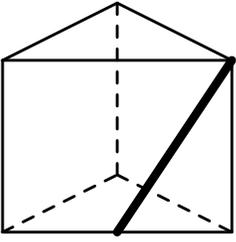 | 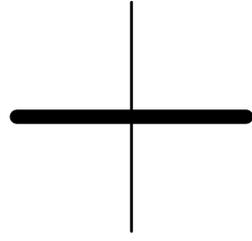<br>1 common double line |
| $S_{\beta_6}$ | 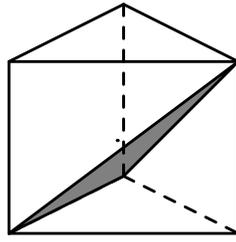 | 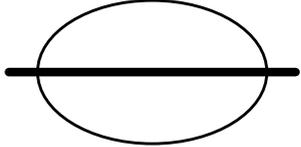<br>2 multiplicity 2 points |
| $S_{\beta_7}$ | 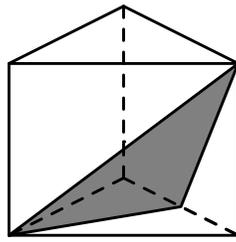 | 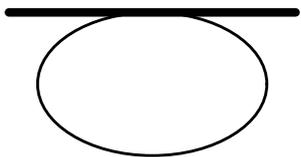<br>The local ring $\cong \bar{k}[\epsilon]/(\epsilon^4)$ |
| $S_{\beta_8}$ | 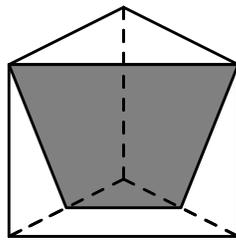 | 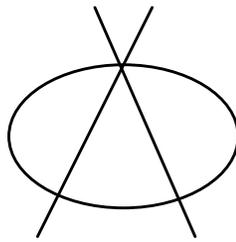<br>1 multiplicity 2 point |

| Strata | Convex hull | Conics |
|---|---|---|
| $S_{\beta_9}$ | 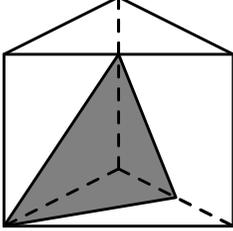 | 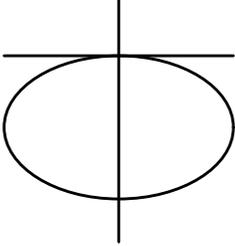<br>1 multiplicity 3 point |
| $S_{\beta_{10}}$ | 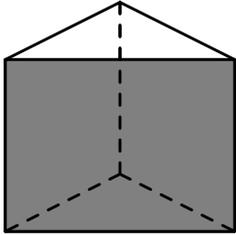 | 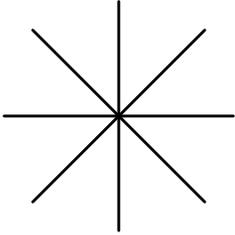<br>The local ring $\cong \bar{k}[\epsilon_1, \epsilon_2]/(\epsilon_1^2, \epsilon_2^2)$ |

Table (2.3)

| $\beta$ | $\|\beta\|^2$ | $Z_\beta$ | $W_\beta$ |
|---|---|---|---|
| $\beta_1$ | $\frac{1}{2}$ | $x_{2,j_1 j_2}$ for $j_1, j_2 = 1, 2, 3$ | — |
| $\beta_2$ | $\frac{7}{6}$ | $x_{2,22}, x_{2,23}, x_{2,33}$ | — |
| $\beta_3$ | $\frac{17}{6}$ | $x_{2,33}$ | — |
| $\beta_4$ | $\frac{1}{6}$ | $x_{1,13}, x_{1,23}, x_{2,13}, x_{2,23}$ | $x_{1,33}, x_{2,33}$ |
| $\beta_5$ | $\frac{11}{12}$ | $x_{1,33}, x_{2,23}$ | $x_{2,33}$ |
| $\beta_6$ | $\frac{1}{24}$ | $x_{1,33}, x_{2,j_1,j_2}$ for $j_1, j_2 = 1, 2, 3$ | $x_{2,13}, x_{2,23}, x_{2,33}$ |
| $\beta_7$ | $\frac{1}{4}$ | $x_{1,33}, x_{2,13}, x_{2,22}$ | $x_{2,23}, x_{2,33}$ |
| $\beta_8$ | $\frac{1}{42}$ | $x_{1,22}, x_{1,23}, x_{1,33}, x_{2,12}, x_{2,13}$ | $x_{2,22}, x_{2,23}, x_{2,33}$ |
| $\beta_9$ | $\frac{1}{10}$ | $x_{1,23}, x_{2,13}, x_{2,22}$ | $x_{1,33}, x_{2,23}, x_{2,33}$ |
| $\beta_{10}$ | $\frac{2}{3}$ | $x_{i,j_1 j_2}$ for $i = 1, 2, j_1, j_2 = 2, 3$ | — |

## §3 An invariant measure on GL(n)

Let $k$ be a number field for the rest of this paper. Let $G = \text{GL}(n)$. We choose an invariant measure on $G$ in this section. We denote the fields of rational, real, complex numbers by $\mathbb{Q}, \mathbb{R}, \mathbb{C}$ respectively. We denote the ring of rational integers by $\mathbb{Z}$. The set of positive real numbers is denoted by $\mathbb{R}_+$. For any ring $R$, $R^\times$ is the set of invertible elements of $R$. Let $k$ be a number field. Let $\mathfrak{M}, \mathfrak{M}_\infty, \mathfrak{M}_\mathbb{R}, \mathfrak{M}_\mathbb{C}, \mathfrak{M}_f$ be the set of all the places, all the infinite, real, imaginary, finite places of $k$ respectively. Let $\mathbb{A}_f$ (resp. $\mathbb{A}_f^\times$) be the restricted product of $k_v$'s (resp. $k_v^\times$'s) over $v \in \mathfrak{M}_f$. Let $k_\infty$ (resp. $k_\infty^\times$) be the product of $k_v$'s (resp. $k_v^\times$'s) over $v \in \mathfrak{M}_\infty$. Then $\mathbb{A} = k_\infty \times \mathbb{A}_f$, $\mathbb{A}^\times = k_\infty^\times \times \mathbb{A}_f^\times$. If $x \in \mathbb{A}$ or $\mathbb{A}^\times$, we denote the finite (resp. infinite) part of $x$ by $x_f$ (resp. $x_\infty$). If $V$ is a vector space over $k$, we define $V_\mathbb{A}, V_\infty, V_f$ similarly. Let $\mathscr{S}(V_\mathbb{A}), \mathscr{S}(V_\infty), \mathscr{S}(V_f)$ be the spaces of Schwartz–Bruhat functions.

For a place $v$, $k_v$ is the completion at $v$. If $v \in \mathfrak{M}_f$, $o_v \subset k_v$ is, by definition, the integer ring of $k_v$. The absolute value of $k_v$ is denoted by $|\ |_v$. For $x \in \mathbb{A}^\times$, we denote the product of $|x|_v$'s over all $v \in \mathfrak{M}_f$ (resp. $v \in \mathfrak{M}_\infty$) by $|x|_f$ (resp. $|x|_\infty$). For $v \in \mathfrak{M}_f$, let $\pi_v$ be the prime element, and $|\pi_v|_v = q_v^{-1}$. Note that if $v$ is imaginary and $|x|$ is the usual absolute value, $|x|_v = |x|^2$.

Let $r_1, r_2$ be the numbers of real and imaginary places respectively. Let $h, R, e$ be the class number, regulator, and the number of units of $k$ respectively. Let $\Delta_k$ be the discriminant of $k$. Let $\mathfrak{C}_k = 2^{r_1}(2\pi)^{r_2} h R e^{-1}$. We choose a Haar measure $dx$ on $\mathbb{A}$ so that $\int_{\mathbb{A}/k} dx = 1$. For any finite place $v$, we choose a Haar measure $dx_v$ on $k_v$ so that $\int_{o_v} dx_v = 1$. We use the ordinary Lebesgue measure $dx_v$ for $v$ real, and $dx_v \wedge d\bar{x}_v$ for $v$ imaginary. Then $dx = |\Delta_k|^{-\frac{1}{2}} \prod_v dx_v$ (see [21, pp. 91]).

For $\lambda \in \mathbb{R}_+$, let $\underline{\lambda}$ be the idele whose component at $v$ is $\lambda^{\frac{1}{[k:\mathbb{Q}]}}$ if $v \in \mathfrak{M}_\infty$ and 1 if $v \in \mathfrak{M}_f$. Clearly, $|\underline{\lambda}| = \lambda$. We identify $\mathbb{R}_+$ with a subgroup of $\mathbb{A}^\times$ by the map $\lambda \to \underline{\lambda}$. Let $\mathbb{A}^1 = \{x \in \mathbb{A}^\times \mid |x| = 1\}$. Then $\mathbb{A}^\times \cong \mathbb{A}^1/k^\times \times \mathbb{R}_+$, and $\mathbb{A}^1/k^\times$ is compact. We choose a Haar measure $d^\times t^1$ on $\mathbb{A}^1$ so that $\int_{\mathbb{A}^1/k^\times} d^\times t^1 = 1$. Using this measure, we choose a Haar measure $d^\times t$ on $\mathbb{A}^\times$ so that

$$\int_{\mathbb{A}^\times} f(t) d^\times t = \int_0^\infty \int_{\mathbb{A}^1} f(\underline{\lambda} t^1) d^\times \lambda d^\times t^1,$$

where $d^\times \lambda = \lambda^{-1} d\lambda$. For any finite place $v$, we choose a Haar measure $d^\times t_v$ on $k_v^\times$ so that $\int_{o_v^\times} d^\times t_v = 1$. Let $d^\times t_v(x) = |x|_v^{-1} dx$ if $v$ is real, and $d^\times t_v(x) = |x|_v^{-1} dx \wedge d\bar{x}$ if $v$ is imaginary. Then $d^\times t = \mathfrak{C}_k^{-1} \prod_v d^\times t_v$ (see [21, pp. 95]). We fix a character $<> = \prod_v <>_v$ of $\mathbb{A}/k$.

Let $\zeta_k(s)$ be the Dedekind zeta function. As in [21], we define

$$Z_k(s) = |\Delta_k|^{\frac{s}{2}} \left(\pi^{-\frac{s}{2}} \Gamma(\frac{s}{2})\right)^{r_1} \left((2\pi)^{-s} \Gamma(s)\right)^{r_2} \zeta_k(s).$$

For a character $\omega$ of $\mathbb{A}^\times/k^\times$, we define $\delta(\omega) = 1$ if $\omega$ is trivial, and $\delta(\omega) = 0$ otherwise.

Let $K = \prod_{v \in \mathfrak{M}} K_v$, where $K_v = O(n)$ if $v$ is real, and $K_v = U(n)$ if $v$ is imaginary, and $K_v = \text{GL}(n, o_v)$ if $v$ is finite. Let $T \subset G$ be the set of diagonal matrices, and $N \subset G$ the set of lower triangular matrices whose diagonal entries are 1. Then $B = TN$ is a Borel subgroup of $G$. We use the notation

$$t = a_n(t_1, \cdots, t_n) = \begin{pmatrix} t_1 & & \\ & \ddots & \\ & & t_n \end{pmatrix},$$

$$n_n(u) = \begin{pmatrix} 1 & & 0 \\ & \ddots & \\ u & & 1 \end{pmatrix}, \; u = (u_{ij})_{i>j} \in \mathbb{A}^{\frac{n(n-1)}{2}},$$

for elements $t \in T_\mathbb{A}, n_n(u) \in N_\mathbb{A}$ respectively. Let

$$T_+ = \{t = a_n(\underline{t}_1, \cdots, \underline{t}_n) \mid t_1, \cdots, t_n > 0\}.$$

The group $G_\mathbb{A}$ has the Iwasawa decomposition $G_\mathbb{A} = KT_\mathbb{A}N_\mathbb{A}$. So any element $g \in G_\mathbb{A}$ can be written as

$$g = k(g)t(g)n_n(u(g)), \; k(g) \in K, \; t(g) \in T_\mathbb{A}, \; n_n(u(g)) \in N_\mathbb{A}.$$

Let $T' = T \cap \mathrm{SL}(n)$, and $T'_+ = T_+ \cap T'_\mathbb{R}$. We consider $\mathfrak{t}$ in §1 for $T'$. We can identify $\mathfrak{t}$ with $\mathfrak{t}_{A_{n-1}}$. For $s \in \mathfrak{t}^*_\mathbb{C}$, we define $t^s, |t|^s_v$ in the usual manner. Let $\rho$ be half the sum of weights of $N$, with respect to the conjugation by elements of $T$. Let $G^0_\mathbb{A} = \{g \in G_\mathbb{A} \mid |\det g| = 1\}$.

Let $du = \prod_{i>j} du_{ij}$. This is an invariant measure on $N_\mathbb{A}$. We choose an invariant measure $dk$ on $K$ so that $\int_K dk = 1$. Let $d^\times t = d^\times t_1 \cdots d^\times t_n$ ($t = a_n(t_1, \cdots, t_n)$). Let $db = t^{-2\rho}d^\times t du$ ($t^{-2\rho} = \prod_{i<j} |t_i^{-1} t_j|$). We choose an invariant measure on $G_\mathbb{A}$ by $dg = t^{-2\rho} dk d^\times t du$. Let $c_n(\lambda) = a_n(\underline{\lambda}, \cdots, \underline{\lambda})$. Let $dg^0$ be the invariant measure on $G^0_\mathbb{A}$ such that for any measurable function on $f(g)$ on $G_\mathbb{A}$,

$$\int_{G_\mathbb{A}} f(g) dg = n \int_0^\infty \int_{G^0_\mathbb{A}} f(c_n(\lambda)g^0) d^\times \lambda dg^0,$$

where $d^\times \lambda = \lambda^{-1} d\lambda$. We define invariant measures on $G_{k_v}, K_v, B_{k_v}, N_{k_v}, T_{k_v}$ similarly, and denote them by $dg_v, dk_v, db_v, du_v, d^\times t_v$ respectively. Then $du = |\Delta_k|^{-\frac{n(n-1)}{4}} \prod_v du_v$, and $d^\times t = \mathfrak{C}_k^{-n} \prod_v d^\times t_v$.

Let

(3.1) $$\mathfrak{V}_n = \frac{Z_k(2) \cdots Z_k(n)}{\mathfrak{R}_k^{n-1}}.$$

We define $\mathfrak{V}_1 = 1$ for convenience. The volume of $\mathrm{GL}(n)^0_\mathbb{A} / \mathrm{GL}(n)_k$ is $\mathfrak{V}_n$.

§4 The zeta function

Let $(G, V, \chi_V)$ be a prehomogeneous vector space. We define the zeta function in a general manner.

Our main interest is a function of the form $Z_L(\Phi, \omega, s)$ in the introduction for irreducible representations, but in order to carry out the induction, we have to consider a slightly more general objects. Let $T_0, \epsilon$ be as in the introduction. For any split torus $T \cong \mathrm{GL}(1)^r$, we define $T_+$ to be the subset of $T_{\mathbb{A}}$ which corresponds to $\mathbb{R}_+^r$ by the above isomorphism. Let $T_+^1 = G_{\mathbb{A}}^1 \cap T_+$.

Let $\omega$ be a character of $G_{\mathbb{A}}^1/G_k$. A principal quasi-character of $G_{\mathbb{A}}^1/G_k$ is a function of the form $|\phi(g)|^s$ where $\phi$ is a rational character of $G$ and $s \in \mathbb{C}$. If $\chi$ is a principal quasi-character of $G_{\mathbb{A}}^1/G_k$, we extend it to $G_{\mathbb{A}}$ so that it it trivial on $T_{0+}$. Let $\chi = (\chi_1, \cdots, \chi_l)$ be principal quasi-characters. We define $\chi(g) = \prod_i \chi_i(g)$. Let $\mathscr{S}(V_{\mathbb{A}})$ be the space of Schwartz–Bruhat functions on $V_{\mathbb{A}}$.

For later purposes, we introduce some notations.

**Definition (4.1)** *Let $L \subset V_k$ be a $G_k$-invariant subset. For $\Phi \in \mathscr{S}(V_{\mathbb{A}})$, we define* $\Theta_L(\Phi, g) = \sum_{x \in L} \Phi(gx)$.

We use the notation $\Theta_{S_\beta}(\Phi, g)$ instead of $\Theta_{S_{\beta k}}(\Phi, g)$ etc.

**Definition (4.2)** *For a $G_k$-invariant subset $L \subset V_k^{\mathrm{ss}}$ and a Schwartz–Bruhat function $\Phi \in \mathscr{S}(V_{\mathbb{A}})$, we define*

$$(1) \qquad Z_L(\Phi, \omega, \chi, s) = \int_{G_{\mathbb{A}}/G_k} \omega(g)\chi(g)|\chi_V(g)|^{\frac{s}{e}} \Theta_L(\Phi, g) dg,$$

$$(2) \qquad Z_{L+}(\Phi, \omega, \chi, s) = \int_{\substack{G_{\mathbb{A}}/G_k \\ |\chi_V(g)| \geq 1}} \omega(g)\chi(g)|\chi_V(g)|^{\frac{s}{e}} \Theta_L(\Phi, g) dg.$$

*if these integtrals converges absolutely.*

For $x \in V_k$, let $G_x^0$ be the connected component of 1 of the stabilizer of $x$. Let $V_k'$ be the set of points $x \in V_k$ such that $G_x^0$ does not have a non-trivial rational character. We expect that the largest $G_k$-invariant set $L$ for which the integrals (1), (2) converges absolutely for $\mathrm{Re}(s) \gg 0$ is $V_k'$.

In [16], F. Sato proved this convergence under the condition that if the stabilizer of $x$ is connected for any $x \in V_k^{\mathrm{ss}}$ and $V_k' = V_k^{\mathrm{ss}}$. F. Sato's paper [16] is written in classical language and the proof in adelic language can be seen in [24]. The above condition applies to many prehomogeneous vector spaces. However, there are several very interesting prehomogeneous vector spaces which do not satisfy the above condition. For those cases, it is useful to use invariant theory as follows.

If $X^*(G)$ is generated by one element, we do not have to consider $\chi$, and therefore we write $Z_L(\Phi, \omega, s)$ etc.

We give an alternative definition of the zeta function when the group is of the form $G/\widetilde{T}$ where $G$ is a product of $\mathrm{GL}(n)$'s, and $\widetilde{T}$ is a split torus contained in the center of $G$.

**Definition (4.3)** *For a connected reductive group $G$, we define $G_{\mathbb{A}}^0 = \{g \in G_{\mathbb{A}} \mid |\chi(g)| = 1 \text{ for all } \chi \in X^*(G)\}$.*

If $G = \mathrm{GL}(n)$, $\mathrm{GL}(n)_{\mathbb{A}}^0 = \{g \in \mathrm{GL}(n)_{\mathbb{A}} \mid |\det g| = 1\}$. Suppose that $G = \mathrm{GL}(n_1) \times \cdots \times \mathrm{GL}(n_f)$. Let $\omega = (\omega_1, \cdots, \omega_f)$ be characters of $\mathbb{A}^\times/k^\times$. We define $\omega(g_1, \cdots, g_f) = \omega_1(\det g_1) \cdots \omega_f(\det g_f)$. Let $\delta_\#(\omega) = \delta(\omega_1) \cdots \delta(\omega_f)$. Let $c_i(\lambda_i)$ be as in §3 for $\lambda_i \in \mathbb{R}_+$. For $\lambda = (\lambda_1, \cdots, \lambda_f) \in \mathbb{R}_+^f$, let $c(\lambda) = (c_1(\lambda_1), \cdots, c_f(\lambda_f))$. Any element of $G_{\mathbb{A}}$ is of the form $c(\lambda)g^0$ where $g^0 \in G_{\mathbb{A}}^0$. We choose a subgroup $\bar{T}_+ \subset T_+^1$ so that $\bar{T}_+ \cong T_+^1/\widetilde{T}_+$. Let $g^1 = \bar{t}g^0$ where $\bar{t} \in \bar{T}$, $g^0 \in G_{\mathbb{A}}^0$. Let $G_{\mathbb{A}}^1 = \bar{T}G_{\mathbb{A}}^0$ and $\widetilde{G}_\mathbb{A} = \mathbb{R}_+ \times G_\mathbb{A}^1$. We identify $T_{0+}$ with $\mathbb{R}_+$. We use the measure $dg^0$ on $G_\mathbb{A}^0$ which is a product of measures we defined in §3. Let $d^\times \bar{t}$ be a Haar measure on $\bar{T}$, and $dg^1 = d^\times \bar{t} dg^0$. Let $\widetilde{g} = (\lambda, g^1)$, where $\lambda \in \mathbb{R}_+$. We define $d\widetilde{g} = d^\times \lambda dg^1$. If we choose $d^\times \bar{t}$ suitably, $Z_L(\Phi, \omega, \chi, s)$ is equal to the following integral

$$(4.4) \qquad \int_{\widetilde{G}_\mathbb{A}/G_k} \lambda^s \omega(g^1) \chi(g^1) \Theta_L(\Phi, \widetilde{g}) d\widetilde{g}.$$

$Z_{L+}(\Phi, \omega, \chi, s)$ has a similar expression also. In particular, if $X^*(G/\widetilde{T})$ is generated by one element, $G_\mathbb{A}^1 = G_\mathbb{A}^0$ and $\bar{T}_+$ is trivial. Therefore, in this case, we use the measure on $G_\mathbb{A}^0$ which we defined in §3. Also we do not have to consider $\chi$. We used this kind of formulation in [26].

We say that $(G, V, \chi_V)$ is of complete type if for any $\chi$, there exists a constant $\sigma(\chi)$ such that $Z_{V_k^{ss}}(\Phi, \omega, \chi, s)$ converges absolutely and locally uniformly for all $\omega, \Phi$, and $\text{Re}(s) > \sigma(\chi)$. Otherwise we say it is of incomplete type. If $(G, V, \chi_V)$ is of complete type, we write $Z_V(\Phi, \omega, \chi, s)$ or simply $Z(\Phi, \omega, \chi, s)$ for $Z_{V_k^{ss}}(\Phi, \omega, \chi, s)$ etc.

Let $G'$ be as in §1, and $T' = T \cap G'$. Let $\mathfrak{t} = X_*(T') \otimes \mathbb{R}$ etc. We choose a coordinate system $x = (x_1, \cdots, x_N)$ whose coordinate vectors are eigenvectors of $T$. Let $\gamma_i \in \mathfrak{t}^*$ be the weights of $x_i$. For $x \in V_k \setminus \{0\}$, we define $I_x = \{i \mid x_i \neq 0\}$. Let $C_x$ be the convex hull of the set $\{\gamma_i \mid i \in I_x\}$.

**Definition (4.5)** *A point $x \in V_k \setminus \{0\}$ is called $k$-stable if for any $g \in G_k$, $C_{gx}$ contains a neighborhood of the origin.*

We denote the set of $k$-stable points by $V_k^s$.

We proved in [26] that $Z_L(\Phi, \omega, \chi, s)$ is well defined for $L = V_k^s$ if $\text{Re}(s)$ is sufficiently large, and $Z_{L+}(\Phi, \omega, \chi, s)$ is an entire function. For the cases where F. Sato's criterion of convergence of the zeta function does not apply, the set $V_k^s$ seems to coincide with $V_k'$. Therefore, except for a few cases, we know the convergence of the zeta function.

Let $G = \text{GL}(n_1) \times \cdots \times \text{GL}(n_f)$, and $(G, V)$ be a prehomogeneous vector space. For $g = (g_1, \cdots, g_f)$, we define ${}^t g = ({}^t g_1, \cdots, {}^t g_f)$. Let $\tau_G$ the longest element of the Weyl group of $G$. Let $V^*$ be the dual space of $V$. The action $f(x) \to f(g^{-1}x)$ for $f \in V^*$ defines a representation of $G$ on $V^*$. The action $x \to {}^t g^{-1} x$ is a representation of $G$ on $V$ and the weights of this representation and the weights of the representation $V^*$ are the same. Therefore, these representations are equivalent. So there exists a bilinear form $[x, y]_V'$ on $V$ such that $[gx, y]_V' = [x, {}^t gy]_V'$. We define $[x, y]_V = [x, \tau_G y]_V'$ and $g^\iota = \tau_G {}^t g^{-1} \tau_G^{-1}$. Then $[gx, y]_V = [x, (g^\iota)^{-1} y]_V$ for all $x, y \in V$ and $g^\iota$ is an involution of $G$.

We define a Fourier transform $\widehat{\Phi}$ of $\Phi \in \mathscr{S}(V_\mathbb{A})$ by the formula

$$(4.6) \qquad \widehat{\Phi}(x) = \int_{V_\mathbb{A}} \Phi(y) <[x, y]_V> dy.$$

The determinant of $\text{GL}(V)$ defines a rational character of $G$ which we call $\bar{\kappa}_V$. Let $\kappa_V$ be the principal quasi-character of $G_\mathbb{A}/G_k$ which is $|\bar{\kappa}_V(g^1)|^{-1}$ for $g^1 \in G_\mathbb{A}^1$, and is trivial on $T_{0+}$. For a $G_k$-invariant subset $L \subset V_\mathbb{A}$, we define

$$(4.7) \qquad J_L(\Phi, \widetilde{g}) = \lambda^{-N} \kappa_V(\widetilde{g}) \sum_{x \in V_k \setminus L} \widehat{\Phi}(\widetilde{g}^\iota x) - \sum_{x \in V_k \setminus L} \Phi(\widetilde{g} x)$$

for $\widetilde{g} = (\lambda, g^1)$. If $L = V_k^{ss}$, we drop $L$, and simply write $J(\Phi, g)$. For $\Phi \in \mathscr{S}(V_\mathbb{A})$ and $\lambda \in \mathbb{R}_+$, we define $\Phi_\lambda(x) = \Phi(\underline{\lambda}x)$. By the Poisson summation formula,

$$(4.8) \quad Z_L(\Phi, \omega, \chi, s) = Z_{L+}(\Phi, \omega, \chi, s) + Z_{L+}(\widehat{\Phi}, \omega^{-1}, \kappa_V^{-1}\chi^{-1}, N - s)$$
$$+ \int_0^1 \int_{G_\mathbb{A}^1/G_k} J_L(\Phi_\lambda, g^1) d^\times \lambda dg^1.$$

For later purposes, we define an operator $M_{V,\omega}$ on $\mathscr{S}(V_\mathbb{A})$ for a character $\omega$ on $\widetilde{G}_\mathbb{A}/G_k$. Consider the standard maximal compact subgroup $K$. We fix a Haar measure $dk$ so that $\int_K dk = 1$.

**Definition (4.9)** *Let $\Phi, \omega$ be as above. We define*

$$M_{V,\omega}\Phi(x) = \int_K \omega(\det k)\Phi(kx)dk.$$

If $\omega$ is trivial, this is the operator to average over $K$.

## §5 The smoothed Eisenstein series

Let $G$ be a connected reductive group not necessarily split over a number field $k$. We fix a minimal parabolic subgroup $P \subset G$. Let $T \subset G$ be a maximal split torus contained in $P$. Let $\rho$ be half the sum of positive weights. Let $d$ be the dimension of $T_+^0$. Let $\mathfrak{t}^0$ be the Lie algebra of $T_+^0$, and $\mathfrak{t}_\mathbb{C}^0 = \mathfrak{t}^0 \otimes \mathbb{C}$. Let $E_P(g, z)$ be the Eisenstein series of $G$ with respect to $P$ where $g \in G_\mathbb{A}^0$, and $z \in \mathfrak{t}_\mathbb{C}^{0*}$.

We fix a Weyl group invariant Hermitian product $(\ ,\ )_0$ on $\mathfrak{t}_\mathbb{C}^{0*}$. Let $\|\ \|_0$ be the norm defined by this product. Let $L(z) = (z, \rho)_0$, and $w_0 = L(\rho)$. We choose an entire function $\psi(z)$ on $\mathfrak{t}_\mathbb{C}^{0*}$ which is rapidly decreasing on any vertical strip.

**Definition (5.1)** *We define a function $\mathscr{E}(g, w, \psi)$ on $G_\mathbb{A}^0/G_k \times \mathbb{C}$ by the following contour integral*

$$\mathscr{E}(g, w, \psi) = \left(\frac{1}{2\pi\sqrt{-1}}\right)^d \int_{\substack{\mathrm{Re}(z)=q \\ \mathrm{Re}(w)>L(q)}} \frac{E_P(g, z)\psi(z)}{w - L(z)} dz,$$

*where we choose $q$ so that $E_P(g, z)$ converges absolutely for $\mathrm{Re}(z) = q$.*

We call $\mathscr{E}(g, w, \psi)$ a smoothed Eisenstein series. If there is no confusion, we drop $\psi$, and write $\mathscr{E}(g, w)$ also.

Let $\alpha_1, \cdots, \alpha_d \in \mathfrak{t}^{0*}$ be simple roots. We fix a small positive number $\eta$. Let $K$ be a special maximal compact subgroup of $G_\mathbb{A}$, so $G_\mathbb{A} = KP_\mathbb{A}$. The group $K$ is contained in $G_\mathbb{A}^0$. Let $T_\eta^0 = \{t \in T_+^0 \mid t^{\alpha_i} \geq \eta \text{ for all } i\}$. We can choose a compact set $\Omega \subset P_\mathbb{A}$ so that $KT_\eta^0\Omega$ surjects to $G_\mathbb{A}^0/G_k$. We can also choose a compact set $\widehat{\Omega} \subset G_\mathbb{A}^0$ so that $\{ktnt^{-1} \mid k \in K, t \in T_\eta^0, n \in \Omega\} \subset \widehat{\Omega}T_\eta^0$.

Let $\mathfrak{t}_+^{0*} = \{\sum_i a_i\alpha_i \mid a_1, \cdots, a_d > 0\}$. If $r_1, r_2 \in \mathfrak{t}^{0*}$, then we use the notation $r_1 > r_2$ if $r_1 - r_2 \in \mathfrak{t}_+^{0*}$. Let $r \in \mathfrak{t}^{0*}$. Let $C(G_\mathbb{A}^0/G_k)$ be the set of continuous functions on $G_\mathbb{A}^0/G_k$. We define

$$(5.2) \quad C(G_\mathbb{A}^0/G_k, r) = \{f \in C(G_\mathbb{A}^0/G_k) \mid \sup_{k \in \widehat{\Omega}, t \in T_\eta} |f(kt)t^r| < \infty\}.$$

If $f \in C(G_{\mathbb{A}}^0/G_k, r)$ for some $r$, we say that $f$ is slowly increasing. A function $f$ is integrable on $G_{\mathbb{A}}^0/G_k$ if $f \in C(G_{\mathbb{A}}^0/G_k, r)$ for some $r > -2\rho$.

**Conjecture (5.3)**

(1) *There exists a constant $C_G$ with the property: for any $\epsilon > 0$, there exists $\delta > 0$ such that if $M > w_0$,*

$$\sup_{\substack{w_0 - \delta \leq \mathrm{Re}(w) \leq M \\ g = kt, \; k \in \widehat{\Omega}, t \in T_\eta}} \left| \mathscr{E}(g, w, \psi) - \frac{C_G \psi(\rho)}{w - w_0} \right| t^r < \infty$$

*for some $r \in \mathfrak{t}^*$ satisfying $\|r\|_0 < \epsilon$.*

(2) *There exist constants $c_1 > 0, c_2$ such that if $M_1 > M_2 > w_0$,*

$$\sup_{\substack{M_2 \leq \mathrm{Re}(w) \leq M_1 \\ g = kt, \; k \in \widehat{\Omega}, t \in T_\eta}} |\mathscr{E}(g, w, \psi)| t^{(c_1 \mathrm{Re}(w) + c_2) \sum_{i=1}^d \alpha_i} < \infty.$$

The statement (1) of the above theorem implies that if $f \in C(G_{\mathbb{A}}^0/G_k, r)$ for some $r > -2\rho$, then the integral

$$\int_{G_{\mathbb{A}}^0/G_k} f(g) \mathscr{E}(g, w, \psi) dg$$

is well defined for $\mathrm{Re}(w) > w_0$, can be continued meromorphically to $\mathrm{Re}(w) > w_0 - \delta$ for some $\delta > 0$, and the principal part around $w = w_0$ is

$$\frac{C_G \psi(\rho)}{w - w_0} \int_{G_{\mathbb{A}}^0/G_k} f(g) dg.$$

The statement (2) of the above theorem implies that if $f$ is a slowly increasing function,

$$\int_{G_{\mathbb{A}}^0/G_k} f(g) \mathscr{E}(g, w, \psi) dg$$

is well defined and becomes a holomorphic function in some right half plane.

If the statement of (5.3) is true for a group $G$, we call it Shintani's lemma for $G$, because we are generalizing Lemma 2.9 [18]. The author proved the following theorem in [26].

**Theorem (5.4)** *The conjecture (5.3) is true if $G$ is a product of $\mathrm{GL}(n)$'s.*

If $G = \mathrm{GL}(n)$, the constant $C_G$ is $\mathfrak{V}_n^{-1}$. The proof of this theorem is fairly technical and involves estimates of Whittaker functions. We briefly explain why we have to use the linear function $L(z)$ by the case of $\mathrm{GL}(n)$. We consider the Borel subgroup $B$.

In this case, the Weyl group is the set of permutations of the set $\{1, \cdots, n\}$. For a permutation $\tau$, we define

$$M_\tau(z) = \prod_{\substack{i > j \\ \tau(i) < \tau(j)}} \frac{Z_k(z_{\tau(i)} - z_{\tau(i)})}{Z_k(z_{\tau(i)} - z_{\tau(i)} + 1)}.$$

The constant term of $E_B(g,z)$ is equal to $\sum_\tau M_\tau(z) t^{\tau z+\rho}$. Let

$$\mathscr{E}_\tau(g,w,\psi) = \left(\frac{1}{2\pi\sqrt{-1}}\right)^{n-1} \int_{\substack{\mathrm{Re}(z)=q \\ \mathrm{Re}(w)>L(q)}} \frac{M_\tau(z) t^{\tau z+\rho}\psi(z)}{w-L(z)} dz,$$

where we use the same $q$ as in (5.1).

Note that if $\tau$ is an element of the Weyl group, $(-\tau^{-1}\rho, -\tau^{-1}\rho)_0 = (\rho,\rho)_0$. Therefore, $-\tau^{-1}\rho$ lies on the sphere of radius $\sqrt{(\rho,\rho)_0}$. It is easy to see that $-\tau^{-1}\rho = \rho$ if and only if $\tau$ is the longest element of the Weyl group. Therefore, $L(-\tau^{-1}\rho) < L(\rho)$ unless $\tau$ is the longest element. If $z = -\tau^{-1}\rho$, then $z_{\tau(i)} - z_{\tau(j)} \geq 1$ for all $(i,j)$ such that $i > j, \tau(i) < \tau(j)$. Therefore, we can choose $q$ in a neighborhood of $-\tau^{-1}\rho$ so that $z_{\tau(i)} - z_{\tau(j)} > 1$ for all such pairs. Then $|\mathscr{E}_\tau(g,w,\psi)| \ll t^{\tau q+\rho}$. The value $\|\tau q + \rho\|_0$ can be arbitrarily small and $L(q) < w_0$, if $\tau$ is not the longest element. This consideration is enough for (1) when $\tau$ is not the longest element. If $\tau$ is the longest element, we consider the residue of $M_\tau(z)$ along $z_1 - z_2 = 1$ then along $z_2 - z_3$ and we continue this process. Since the eventual residue is the constant function $\mathfrak{V}_n^{-1}$ as a function of $g$, we get (1). This argument for constant terms is fairly easy and is seen in Langlands' work [11]. However, in [26], we spent most of our labor to the estimates of the non-constant terms. The statement (2) is slightly trickier, and the reader should see [26] for the details of the proof of Theorem (5.4).

## §6 Adjusting terms

In this section, we illustrate the necessity of adjusting terms by the example of the space of binary quadratic forms. Let $G = \mathrm{GL}(1) \times \mathrm{GL}(2)$, $V = \mathrm{Sym}^2 k^2$. Let $\widetilde{T} \subset G$ be the kernel of the homomorphism $G \to \mathrm{GL}(V)$. Let $\chi_V(t,g) = t \det g$ for $(t,g) \in G$. Then $\chi_V$ is a rational character of $G/\widetilde{T}$ and is indivisible as a character of $G$. In this case, $V_k^{\mathrm{s}}$ is the set of forms without rational factors.

Let $\omega$ be as in §3. We do not have to consider $\chi$ in (4.1). We use the alternative definition (4.3) of the zeta function for $L = V_k^{\mathrm{s}}$ i.e.,

(6.1) $$Z_V(\Phi,\omega,s) = \int_{\mathbb{R}_+ \times G_\mathbb{A}^0/G_k} \lambda^s \omega(g^0) \Theta_{V_k^{\mathrm{s}}}(\Phi, \underline{\lambda} g^0) d^\times \lambda \, dg^0 \text{ etc.}$$

**Definition (6.2)** *For $u \in \mathbb{A}$, we define*

$$\alpha(u) = \prod_{v \in \mathfrak{M}} \alpha_v(u_v), \text{ where } \alpha_v(u_v) = \begin{cases} \sup(1, |u_v|_v)^{-1} & \text{for } v \in \mathfrak{M}_f, \\ (1+|u_v|_v^2)^{-\frac{1}{2}} & \text{for } v \in \mathfrak{M}_\mathbb{R}, \\ (1+|u_v|_v)^{-1} & \text{for } v \in \mathfrak{M}_\mathbb{C}. \end{cases}$$

**Definition (6.3)** *Let $v$ be a place of $k$. Let $\Psi, \Psi_v$ be Schwartz–Bruhat functions on $\mathbb{A}^2, k_v^2$ respectively. Let $s, s_1$ be complex variables, and $\omega$ a character of $\mathbb{A}^\times/k^\times$.*

*We define*

(1) $$T_{V,v}(\Psi_v, \omega, s, s_1) = \int_{k_v^\times \times k_v} \omega(t_v)|t_v|_v^s \alpha_v(u_v)^{s_1} \Psi_v(t_v, t_v u_v) d^\times t_v du_v,$$

(2) $$T_V(\Psi, \omega, s, s_1) = \int_{\mathbb{A}^\times \times \mathbb{A}} \omega(t)|t|^s \alpha(u)^{s_1} \Psi(t, tu) d^\times t du,$$

(3) $$T_{V+}(\Psi, \omega, s, s_1) = \int_{\substack{\mathbb{A}^\times \times \mathbb{A} \\ |t| \geq 1}} \omega(t)|t|^s \alpha(u)^{s_1} \Psi(t, tu) d^\times t du,$$

(4) $$T_V^1(\Psi, \omega, s_1) = \int_{\mathbb{A}^1 \times \mathbb{A}} \omega(t^1) \alpha(u)^{s_1} \Psi(t^1, t^1 u) d^\times t^1 du,$$

(5) $$T_{V,v}(\Psi_v, \omega, s, s_1) = \frac{d}{ds_1}\bigg|_{s_1=0} T_V(\Psi, \omega, s, s_1).$$

Let $\widetilde{R}_{V,0}\Phi \in \mathscr{S}(\mathbb{A}^2)$ be the restriction of $\Phi$ to the set of elements of the form $(0, x_1, x_2)$

**Definition (6.4)**

$$Z_{V,\mathrm{ad}}(\Phi, \omega, s) = Z(\Phi, \omega, s) - \frac{\delta(\omega_1 \omega_2^{-1})}{2} T_V(\widetilde{R}_{V,0} M_{V,w}\Phi, \omega_1, s)$$

If $\Psi = \otimes \Psi_v$, then $T_V(\Psi, \omega, s, s_1) = |\Delta_k|^{-\frac{1}{2}} \mathfrak{C}_k^{-1} \prod_v T_v(\Psi_v, \omega, s, s_1)$. Suppose that $\Psi_v$ is the characteristic function of $o_v^2$, and $\omega$ is trivial. Then easy computations (see (2.10) [25]) show that

$$T_v(\Psi_v, \omega, s, s_1) = \frac{1 - q_v^{-(s+s_1)}}{(1 - q_v^{-s})(1 - q_v^{-(s+s_1-1)})}.$$

In general, $T_v(\Psi_v, \omega, s, s_1)$ can be continued meromorphically everywhere, and is holomorphic for $\mathrm{Re}(s) + \mathrm{Re}(s_1) > 1, \mathrm{Re}(s) > 0$. Therefore, $T_V(\Psi, \omega, s, s_1)$ can be continued meromorphically everywhere, and is holomorphic for $\mathrm{Re}(s) + \mathrm{Re}(s_1) > 2, \mathrm{Re}(s) > 1$. Also the integrals defining $T_{V+}(\Psi, \omega, s, s_1)$, $T_V^1(\Psi, \omega, s_1)$ converge absolutely and locally uniformly for all $s, s_1$. We call $Z_{V,\mathrm{ad}}(\Phi, \omega, s)$ the adjusted zeta function for $V$, and $\frac{\delta(\omega_1 \omega_2^{-1})}{2} T_V(\widetilde{R}_{V,0} M_{V,w}\Phi, \omega_1, s)$ the adjusting term.

$Z_{V,\mathrm{ad}}(\Phi, \omega, s)$ has finitely many poles and satisfies a functional equation

$$Z_{V,\mathrm{ad}}(\Phi, \omega, s) = Z_{V,\mathrm{ad}}(\widehat{\Phi}, \omega^{-1}, 3 - s).$$

We will state the principal part of $Z_{V,\mathrm{ad}}(\Phi, \omega, s)$ in §8. Also the author proved that generalized special values of $Z_{V,\mathrm{ad}}(\Phi, \omega, s)$ contribute to the poles of zeta functions in [26].

First we briefly explain how the adjusting term arises. Let $H_V$ be the subgroup of $G$ generated by $T$ and the element $(1, \tau_{\mathrm{GL}(2)})$.

**Proposition (6.5)** *Let*

$$X_{H_V} = \{(t, ka(t_1, t_2)n(u)) \in G_\mathbb{A}^0 \mid k \in K, |t_1 t_2^{-1}| \geq \sqrt{\alpha(t_1^{-1} t_2 u)}\}.$$

*Then for any measurable function $f(g)$ on $G^0_{\mathbb{A}}/H_{Vk}$,*

$$\int_{G^0_{\mathbb{A}}/H_{Vk}} f(g)dg = \int_{X_{H_V}/T_k} f(g)dg,$$

*if the right hand side converges absolutely.*

Let $Z'^{\text{ss}}_0$ be the set of points of the form $(0, x_1, 0)$ such that $x_1 \neq 0$. Then $G_k Z'^{\text{ss}}_{0,k} = G_k \times_{H_{Vk}} Z'^{\text{ss}}_{0,k}$. The reason why we need the adjusting term in this case is that we have to consider the fundamenental domain $X_{H_V}$ for the group $H_V$ in order to handle the set $V^{\text{ss}}_k \setminus V^{\text{s}}_k$.

Suppose that $(G, V, \chi_V)$ is a prehomogeneous vector space of incomplete type. Then if $L$ is the maximal subset such that (4.1) converges absolutely and locally uniformly for $\text{Re}(s) \gg 0$ (assuming such $L$ exists), the author expects that $V^{\text{ss}}_k \setminus L$ decomposes into finite $G_k$-invariant subsets $S_{st,1}, \cdots, S_{st,i_V}$ where $S_{st,i} \cong G_k \times_{H_{V,ik}} Z'^{\text{ss}}_{0,ik}$, $H_{V,i} \subset G$ is a (not necessarily connected) reductive subgroup of $G$, $Z'_{0,i} \subset V$ is an $H_{V,i}$-invariant subspace, and $Z'^{\text{ss}}_{0,ik}$ is the set of semi-stable points with respect to $H_{V,i} \cap G'$ for all $i$. Let $H^0_{V,i}$ be the connected component of $H_{V,i}$ which contains the unit element. Suppose that for each $i$, there exist functions $\alpha_{V,i1}, \cdots, \alpha_{V,ij_{V,i}}$ and that if

$$X_{H_{V,i}} = \{g^1 \in G^1_{\mathbb{A}} \mid \alpha_{V,i1}(g^1) \geq 1, \cdots, \alpha_{V,ij_{V,i}}(g^1) \geq 1\},$$

then $X_{H_{V,i}}$ is invariant under the right action of $H^0_{V,ik}$, and

$$\int_{G^1_{\mathbb{A}}/H_{V,ik}} f(g^1)dg^1 = \int_{X_{H_{V,i}}/H^0_{V,ik}} f(g^1)dg^1.$$

Suppose that these conditions are satisfied.

**Definition (6.6)** *Let $s_i = (s_{i,1}, \cdots, s_{i,j_{V,i}}) \in \mathbb{C}^{j_{V,i}}$. We define*

(1) $\quad T_{V,i}(\Phi, \omega, \chi, s, s_i)$

$$= \int_{\mathbb{R}_+ \times G^1_{\mathbb{A}}/H_{V,ik}} \lambda^s \omega(g^1) \chi(g^1) \prod_j \alpha_{ij}(g^1)^{s_j} \Theta_{S_{st,i}}(\Phi, \underline{\lambda}g^1) d^{\times}\lambda dg^1,$$

(2) $\quad T_{V,i+}(\Phi, \omega, \chi, s, s_i)$

$$= \int_{\substack{\mathbb{R}_+ \times G^1_{\mathbb{A}}/H_{V,ik} \\ \lambda \geq 1}} \lambda^s \omega(g^1) \chi(g^1) \prod_j \alpha_{ij}(g^1)^{s_j} \Theta_{S_{st,i}}(\Phi, \underline{\lambda}g^1) d^{\times}\lambda dg^1,$$

(3) $\quad T^1_{V,i}(\Phi, \omega, \chi, s_i)$

$$= \int_{G^1_{\mathbb{A}}/H_{V,ik}} \omega(g^1) \chi(g^1) \prod_j \alpha_{ij}(g^1)^{s_j} \Theta_{S_{st,i}}(\Phi, g^1) dg^1,$$

(4) $\quad T_{V,i}(\Phi, \omega, \chi, s)$

$$= \frac{d}{ds_1}\bigg|_{s_1=0} \cdots \frac{d}{ds_{j_{V,i}}}\bigg|_{s_{j_{V,i}}=0} T_{V,i}(\Phi, \omega, \chi, s, s_1, \cdots, s_{j_{V,i}}) \text{ etc.}$$

*if these functions are well defined.*

We expect $T_{V,i+}(\Phi,\omega,\chi,s,s_i), T^1_{V,i}(\Phi,\omega,\chi,s_i)$ to be entire and $T_{V,i}(\Phi,\omega,\chi,s,s_i)$ is well defined if $\mathrm{Re}(s) \gg 0$ depending on $s_i$.

In the case $V = \mathrm{Sym}^2 k^2$, our adjusting term appeared in this way. Also we proved in [25] that $T_{V+}(\Psi,\omega,s)$ is an entire function, and $T_V(\Psi,\omega,s)$ is holomorphic for $\mathrm{Re}(s) > 2$. Therefore, $T_V(\Psi,\omega,s)$ does not effect the right most pole. The author expects similar properties for the general adjusting terms (see [26] for example).

### §7 The general program

In this section, we describe our program to compute the principal parts of zeta functions for prehomogeneous vector spaces when $G$ is a product of $\mathrm{GL}(n)$'s.

**Definition (7.1)** *Let $f_1(w), f_2(w)$ be meromorphic functions on a domain of the form $\{w \in \mathbb{C} \mid \mathrm{Re}(w) > A\}$. We use the notation $f_1 \sim f_2$ if the following two conditions are satisfied.*
*(1) There exists a constant $A' < w_0$ such that $f_1 - f_2$ can be continued meromorphically to the domain $\{w \in \mathbb{C} \mid \mathrm{Re}(w) > A'\}$.*
*(2) $f_1(w) - f_2(w)$ is holomorphic around $w = w_0$.*

Consider the identification $G^1_{\mathbb{A}}/G_k \cong \bar{T} \times G^0_{\mathbb{A}}/G_k$ in §4. Suppose that

$$|J_L(\Phi, \bar{t}kt^0)| \ll h_1(\bar{t})|t^0|^r,$$

where $h_1(\bar{t})$ is an integrable function on $\bar{T}$, $kt^0 \in \widehat{\Omega}T_\eta$, and that $r \in \mathfrak{t}^{0*}$ satisfies the condition $r > -2\rho$. Let $dg^1 = d^{\times}\bar{t}dg^0$. By (5.4),

$$\int_{\bar{T} \times G^0_{\mathbb{A}}/G_k} \omega(g^0)\chi(\bar{t})J_L(\Phi, \bar{t}g^0)\mathscr{E}(g^0, w, \psi)dg^0$$
$$\sim \frac{C_G \psi(\rho)}{w - w_0} \int_{\bar{T} \times G^0_{\mathbb{A}}/G_k} \omega(g^0)\chi(\bar{t})J_L(\Phi, \bar{t}g^0)dg^0.$$

Suppose that (6.6) is well defined. Then we can state our program as follows.

**Program (7.2)** *Find constants $C_{V,\mathrm{st},1}, \cdots, C_{V,\mathrm{st},i_V}$, points $p_1(\chi), \cdots, p_{d_V}(\chi) \in \mathbb{C}$ and distributions $a(\Phi,\omega,\chi), a_{ij}(\Phi,\omega,\chi)$ for $i = 1, \cdots, d_V, j = 1, \cdots, d_{V,i}$ with the following properties*

(1) $\displaystyle\int_{G^1_{\mathbb{A}}/G_k} \omega(g^0)\chi(\bar{t})J_L(\Phi, g^1)\mathscr{E}(g^1, w, \psi)dg^1$

$\displaystyle\sim \frac{C_G \psi(\rho)}{w - w_0} \sum_{i=1}^{i_V} C_{V,\mathrm{st},i}(T^1_{V,i}(\widehat{\Phi}, \omega^{-1}) - T^1_{V,i}(\Phi,\omega)) + \frac{\psi(\rho)}{w - w_0} a(\Phi,\omega,\chi),$

(2) $\displaystyle a(\Phi_\lambda, \omega, \chi) = \sum_{i=1}^{d_V} \sum_{j=1}^{d_{V,i}} \lambda^{-p_i(\chi)}(\log \lambda)^j a_{ij}(\Phi,\omega,\chi),$

*where $\lambda \in \mathbb{R}_+$.*

The following equations are easy to verify.

$$\int_0^1 \lambda^s T_{V,i}^1(\widehat{\Phi}_\lambda, \omega^{-1}, \kappa_V^{-1}\chi^{-1}) d^\times \lambda = T_{V,i+}(\widehat{\Phi}, \omega^{-1}, \kappa_V^{-1}\chi^{-1}, N-s),$$

$$\int_0^1 \lambda^s T_{V,i}^1(\Phi_\lambda, \omega, \chi) d^\times \lambda = T_{V,i}(\Phi, \omega, \chi, s) - T_{V,i+}(\Phi, \omega, \chi, s).$$

Also
$$\int_0^1 \lambda^{s-p_i(\chi)}(\log \lambda)^j d^\times \lambda = (-1)^j j! (s - p_i(\chi))^{-j-1}.$$

We define

$$(7.3) \qquad Z_{L,\mathrm{ad}}(\Phi, \omega, \chi, s) = Z_L(\Phi, \omega, \chi, s) + \sum_{i=1}^{d_V} C_{st,i} T_{V,i}(\Phi, \omega, \chi, s)$$

and call $Z_{L,\mathrm{ad}}(\Phi, \omega, \chi, s)$ the adjusted zeta function.

Then if (7.2) is carried out, we get the following formula for $Z_{L,\mathrm{ad}}(\Phi, \omega, \chi, s)$

$$(7.4) \quad Z_{L,\mathrm{ad}}(\Phi, \omega, \chi, s)$$
$$= Z_{L+}(\Phi, \omega, \chi, s) + Z_{L+}(\widehat{\Phi}, \omega^{-1}, \kappa_V \chi^{-1}, N-s)$$
$$+ \sum_{i=1}^{i_V} C_{st,i}(T_{V,i+}(\Phi, \omega, \chi, s) + T_{V,i+}(\widehat{\Phi}, \omega^{-1}, \kappa_V \chi^{-1}, N-s))$$
$$+ C_G^{-1} \sum_{i=1}^{d_V} \sum_{j=1}^{d_{V,i}} (-1)^j j! \frac{a_{ij}(\Phi, \omega, \chi)}{(s - p_i(\chi))^{j+1}}.$$

We expect that the first four terms to be entire and that the adjusting terms are holomorphic at $s = \max p_i(\chi)$.

The author carried out this program for $\mathrm{Sym}^2 k^n$ and $\mathrm{Sym}^2 k^3 \otimes k^2$. We describe our results in §8. Consider the second case. The stratum $S_{\beta_8}$ in §2 corresponds to the representation (2) in §2. The zeta function for this case was handled by F. Sato (see [14]) in a slightly different formulation. The author handled two more reducible prehomogeneous vector spaces which correspond to $\beta_6, \beta_8$ in §2 from our viewpoint. These are reducible representations with two irreducible components. The reader should see [26] for the details of these cases. In the future, we have to generalize the method so that we can handle representations with more than two irreducible factors.

For the rest of this section, we briefly discuss how we proceed, and point out technical difficulties. We consider an irreducible prehomogeneous vector space of complete type $(G, V, \chi_V)$ for simplicity. So $G_{\mathbb{A}}^1 = G_{\mathbb{A}}^0$, and we do not have to consider $\chi$.

By the Morse stratification,

$$J(\Phi, g^1) = \sum_\beta \left( \Theta_{S_\beta}(\widehat{\Phi}, (g^1)^\iota) - \Theta_{S_\beta}(\Phi, g^1) \right) + \widehat{\Phi}(0) - \Phi(0).$$

($\kappa_V(g^1) = 1$.)

Therefore,

$$\int_{G^1_{\mathbb{A}}/G_k} \omega(g^1) J(\Phi, g^1) \mathscr{E}(g^1, w, \psi) dg^1$$

$$= \sum_\beta \int_{G^1_{\mathbb{A}}/G_k} \omega(g^1) \Theta_{S_\beta}(\widehat{\Phi}, (g^1)^\iota) \mathscr{E}(g^1, w, \psi) dg^1$$

$$- \sum_\beta \int_{G^1_{\mathbb{A}}/G_k} \omega(g^1) \Theta_{S_\beta}(\Phi, g^1) \mathscr{E}(g^1, w, \psi) dg^1$$

$$+ \frac{\psi(\rho)}{w - w_0}(\widehat{\Phi}(0) - \Phi(0)).$$

In this way, we can at least start considering the contributions from unstable strata indirectly by the use of the smoothed Eisenstein series. Let

$$\Xi_\beta(\Phi, \omega, w) = \int_{G^1_{\mathbb{A}}/G_k} \omega(g^1) \Theta_{S_\beta}(\Phi, g^1) \mathscr{E}(g^1, w, \psi) dg^1.$$

First, it can happen that $\Xi_\beta(\Phi, \omega, w) \sim 0$ for some strata $\beta$. For $\mathrm{Sym}^2 k^3 \otimes k^2$, this is the case for $\beta = \beta_2, \beta_3, \beta_5$. Therefore, not all the strata contribute to the poles, and these strata have to be removed.

If $\Phi = M_{V,\omega}\Phi$,

$$\Xi_\beta(\Phi, \omega, w) = \int_{P^1_{\beta\mathbb{A}}/P_{\beta k}} \omega(g^1_\beta) \Theta_{Y^{ss}_\beta}(\Phi, g^1_\beta n_\beta(u_\beta)) \mathscr{E}(g^1_\beta n_\beta(u_\beta), w, \psi)(g^1_\beta)^{-2\rho_\beta} dg^1_\beta du_\beta,$$

where $P^1_{\beta\mathbb{A}} = P_{\beta\mathbb{A}} \cap G^1_{\mathbb{A}}$, $\rho_\beta$ is half the sum of positive weights of $P_\beta$, and $dg^1_\beta, du_\beta$ are invariant measures on $M_{\beta\mathbb{A}} \cap G^1_{\mathbb{A}}, U_{\beta\mathbb{A}}$. We want to compare this integral with the integration on $M_{\beta\mathbb{A}} \cap G^1_{\mathbb{A}}/M_{\beta k}$ by considering the constant terms of both $\Theta_{Y^{ss}_\beta}(\Phi, g^1_\beta n_\beta(u_\beta))$, and $\mathscr{E}(g^1_\beta n_\beta(u_\beta), w, \psi)$ with respect to $U_\beta$. Then the resulting integral resembles the zeta function for $(M_\beta, Z_\beta)$ except it still has the Eisenstein series in the integrand. However, the estimate of the non-constant terms is a formidable task. In the case $\mathrm{Sym}^2 k^3 \otimes k^2$, this was done by considering further stratifications of the representations $Z_\beta$ using the notion of $\beta$-sequences ([10, pp. 73]).

Also it is too optimistic to expect that for all $\beta$,

$$\int_{G^1_{\mathbb{A}}/G_k} \omega(g^1) \Theta_{S_\beta}(\Phi, g^1) \mathscr{E}(g^1, w, \psi) dg^1 \sim \frac{\psi(\rho)}{w - w_0} a_\beta(\Phi, \omega)$$

for some $a_\beta(\Phi, \omega)$, satisfying a similar property as in (7.2) (2). This is because there are interactions (both constant terms and non-constant terms) between strata. For example, for $\mathrm{Sym}^2 k^3 \otimes k^2$, some cancellations have to be established between the strata $S_{\beta_6}, S_{\beta_8}, S_{\beta_9}$. Those cancellations can be established only if we consider the further Morse stratifications of the representations $Z_{\beta_8}$ etc. The reader should see [26] for the details.

## §8 Examples of principal part formulas

In this section, we consider some examples of prehomogeneous vector spaces for which the principal parts of the poles of the zeta function were determined. In order to state them, we have to introduce some notation.

Let $\Phi \in \mathscr{S}(V_\mathbb{A})$, and $S_\beta = G \times_{P_\beta} Y_\beta^{\mathrm{ss}}$ be a stratum. We define $R_\beta \Phi \in \mathscr{S}(Z_{\beta\mathbb{A}})$ by the formula

$$R_\beta \Phi(x) = \int_{W_{\beta\mathbb{A}}/W_{\beta k}} \Phi(x,y) dy$$

for $x \in Z_{\beta\mathbb{A}}$, where $dy$ is the measure on $W_{\beta\mathbb{A}}$ such that $\int_{W_{\beta\mathbb{A}}/W_{\beta k}} dy = 1$.

Let $\Psi \in \mathscr{S}(\mathbb{A}^i)$. Let $\omega_1, \cdots, \omega_i$ be characters of $\mathbb{A}^\times/k^\times$, and $\omega = (\omega_1, \cdots, \omega_i)$. For $t = (t_1, \cdots, t_i) \in \mathbb{A}^i$, we define $\omega(t) = \prod_j \omega_j(t_j)$. We consider $\omega$ as a character of $(\mathbb{A}^\times/k^\times)^i$. Let $\Theta_i(\Psi, t) = \sum_{x \in (k^\times)^i} \Psi(tx)$. For $s = (s_1, \cdots, s_i) \in \mathbb{C}^i$, and $t$ as above, we define $|t|^s$ in an obvious manner.

**Definition (8.1)**

$$\Sigma_i(\Psi, \omega, s) = \int_{(\mathbb{A}^\times/k^\times)^i} |t|^s \omega(t) \Theta_i(\Psi, t) d^\times t.$$

If $\omega$ is trivial, we drop $\omega$, and use the notation $\Sigma_i(\Psi, s)$.

**Definition (8.2)** *For $j = (j_1, \cdots, j_i) \in \mathbb{Z}^i$, $s_0 = (s_{0,1}, \cdots, s_{0,i}) \in \mathbb{C}^i$, we use the notation $\Sigma_{i,(j_1, \cdots, j_i)}(\Psi, \omega, s_0)$ for the coefficient of $\prod_l (s_l - s_{0,l})^{j_l}$ in the Laurent expansion of $\Sigma_i(\Psi, \omega, s)$ at $s = s_0$. We also drop $\omega$ in this notation when $\omega$ is trivial.*

If $\Sigma_i(\Psi, \omega, s)$ is holomorphic at $s = s_0$, $\Sigma_{i,(0,\cdots,0)}(\Psi, \omega, s_0) = \Sigma_i(\Psi, \omega, s_0)$.

(1) The case $G = \mathrm{GL}(1) \times \mathrm{GL}(2)$, $V = \mathrm{Sym}^3 k^2$

We consider the situation in the second example in §2. We define the zeta function by (4.4) (see the remark after (4.4)) for $L = V_k^{\mathrm{ss}}$. The following theorem was proved by Shintani [18], and the adelic version can be found in [22].

**Theorem (8.3)** *Suppose that $\Phi = M_{V,\omega}\Phi$. Then*

$$Z_V(\Phi, \omega, s) = Z_{V+}(\Phi, \omega, s) + Z_{V+}(\widehat{\Phi}, \omega^{-1}, 4-s) + \mathfrak{V}_2 \delta_\#(\omega) \left( \frac{\widehat{\Phi}(0)}{s-4} - \frac{\Phi(0)}{s} \right)$$

$$+ \delta_\#(\omega) \left( \frac{\Sigma_1(R_{\beta_1}\widehat{\Phi}, -1)}{s-4} - \frac{\Sigma_1(R_{\beta_1}\Phi, -1)}{s} \right)$$

$$+ \delta(\omega_1^3)\delta(\omega_2) \left( \frac{\Sigma_1(R_{\beta_2}\widehat{\Phi}, \omega_1^{-1}, \frac{2}{3})}{3s-10} - \frac{\Sigma_1(R_{\beta_2}\Phi, \omega_1, \frac{2}{3})}{3s-2} \right).$$

(2) The case $G = \mathrm{GL}(1) \times \mathrm{GL}(2)$, $V = \mathrm{Sym}^2 k^2$

In this case, there is precisely one unstable stratum $S_{\beta_1}$, where $\beta_1 = 2$ using a similar identification $\mathbb{R} \cong \mathfrak{t}$ as in §2 (1). The following theorem was proved by Shintani [19], and the adelic version can be found in [25].

**Theorem (8.4)** *Suppose that* $\Phi = M_{V,\omega}\Phi$. *Then*
$$Z_{V,\mathrm{ad}}(\Phi,\omega,s) = Z_{V+}(\Phi,\omega,s) + Z_{V+}(\widehat{\Phi},\omega^{-1},3-s)$$
$$- \frac{1}{2}T_{V+}(\widetilde{R}_{V,0}\Phi,\omega,s) - \frac{1}{2}T_{V+}(\widetilde{R}_{V,0}\widehat{\Phi},\omega^{-1},3-s)$$
$$+ \mathfrak{V}_2\delta_\#(\omega)\left(\frac{\widehat{\Phi}(0)}{s-3} - \frac{\Phi(0)}{s}\right)$$
$$+ \frac{\delta(\omega_1^2)\delta(\omega_2)}{2}\left(-\frac{\Sigma_{1,(-1)}(R_{\beta_1}\widehat{\Phi},\omega_1^{-1},1)}{(s-2)^2} + \frac{\Sigma_{1,(0)}(R_{\beta_1}\widehat{\Phi},\omega_1^{-1},1)}{s-2}\right)$$
$$- \frac{\delta(\omega_1^2)\delta(\omega_2)}{2}\left(\frac{\Sigma_{1,(-1)}(R_{\beta_1}\Phi,\omega_1,1)}{(s-1)^2} + \frac{\Sigma_{1,(0)}(R_{\beta_1}\Phi,\omega_1,1)}{s-1}\right).$$

(3) The case $G = \mathrm{GL}(1) \times \mathrm{GL}(n)$, $V = \mathrm{Sym}^2 k^n$

The case $n = 2$ is the previous case. Let $\widetilde{T}$ be as before, and we consider $(G/\widetilde{T}, V)$. If $n \geq 3$, it is of complete type and $X^*(G/\widetilde{T})$ is generated by one element. So we consider the situation in (4.4). We define
$$Z_{V_n,(j)}(\Phi,\omega,s_0),\ Z_{V_2,\mathrm{ad},(j)}(\Phi,\omega,s_0)$$
similarly as in (8.2).

Let $S_{\beta_i}$ be the set of rank $n-i$ forms. Then $S_{\beta_1},\cdots,S_{\beta_{n-1}}$ are unstable strata. For all $i$, $Y_{\beta_i} = Z_{\beta_i} \cong V_{n-i}$. Let $Z'_{\beta_{n-2},0} \subset Z_{\beta_{n-2}}$ be the set of points of the form $\{(0,x_1,0)\}$ by the above identification. For $\Phi \in \mathscr{S}(V_\mathbb{A})$ and $x_1 \in Z'_{\beta_{n-2}}$, let
$$R'_{\beta_{n-2},0}\Phi(x_1) = \int_\mathbb{A} R_{\beta_{n-2}}\Phi(0,x_1,x_2)dx_2.$$

Shintani proved a partial result, and the author carried out the complete calculation of the principal parts in this case (see [19], [25]).

If $n = 3$, we define
$$F_3(\Phi,\omega,s) = \mathfrak{V}_3\delta_\#(\omega)\frac{\Phi(0)}{s} + \mathfrak{V}_2\delta(\omega_2)\frac{\Sigma_1(R_{\beta_2}\Phi,\omega_1,\frac{3}{2})}{2s-3}$$
$$+ \delta(\omega_2)\frac{Z_{V_2,\mathrm{ad},(0)}(R_{\beta_1}\Phi,(\omega_1,1),3)}{s-3}$$
$$+ \delta_\#(\omega)\left(\frac{3\Sigma_1(R'_{\beta_1,0}\Phi,2)}{2(s-3)^2} + \frac{\Sigma_{1,(1)}(R'_{\beta_1,0}\Phi,2)}{2(s-3)}\right).$$

If $n \geq 4$, we define
$$F_n(\Phi,\omega,s) = \mathfrak{V}_n\delta_\#(\omega)\frac{\Phi(0)}{s} + \sum_{i\neq 1,n-2}(n-i)\mathfrak{V}_i\delta(\omega_2)\frac{Z_{V_{n-i}}(R_{\beta_i}\Phi,(\omega_1,1),\frac{n(n-i)}{2})}{2s-n(n-i)}$$
$$+ \mathfrak{V}_{n-2}\delta(\omega_2)\frac{Z_{V_2,\mathrm{ad}}(R_{\beta_{n-2}}\Phi,(\omega_1,1),n)}{s-n}$$
$$+ (n-1)\delta(\omega_2)\frac{Z_{V_{n-1},(0)}(R_{\beta_1}\Phi,(\omega_1,1),\frac{n(n-1)}{2})}{2s-n(n-1)}$$
$$+ \mathfrak{V}_{n-2}\delta_\#(\omega)\left(\frac{n\Sigma_1(R'_{\beta_{n-2},0}\Phi,n-1)}{2(s-n)^2} + \frac{\Sigma_{1,(1)}(R'_{\beta_{n-2},0}\Phi,n-1)}{2(s-n)}\right).$$

**Theorem (8.5)** *Let $n \geq 3$. Suppose that $M_{V_n,\omega}\Phi = \Phi$. Then*

$$Z_{V_n}(\Phi,\omega,s) = Z_{V_n+}(\Phi,\omega,s) + Z_{V_n+}(\widehat{\Phi},\omega^{-1},\frac{n(n+1)}{2} - s)$$
$$- F_n(\widehat{\Phi},\omega^{-1},\frac{n(n+1)}{2} - s) - F_n(\Phi,\omega,s).$$

Suppose that $\Phi = \otimes_v \Phi_v$, and $\Phi_v$ has a compact support contained in $V_{k_v}^{\mathrm{ss}}$ for a place $v \in \mathfrak{M}_\infty$. Then it turns out that $\Sigma_1(R'_{\beta_{n-2},0}\Phi, n-1) = \Sigma_1(R'_{\beta_{n-2},0}\widehat{\Phi}, n-1) = 0$. Therefore, the poles of the associated Dirichlet series are all simple.

(4) The case $G = \mathrm{GL}(3) \times \mathrm{GL}(2)$, $V = \mathrm{Sym}^2 k^3 \otimes k^2$

This is the third example in §2. Let $\omega = (\omega_1,\omega_2)$ be as before. In order to handle the zeta function in this case, we need to know the pole structures of zeta functions for unstable strata. The strata $S_{\beta_2}, S_{\beta_3}, S_{\beta_5}$ do not contribute to the poles. The stratum $S_{\beta_1}$ corresponds to $V_3$ in (3). Let $\Sigma_{\beta_1}(\Psi,\omega',s) = Z_{V_3}(\Psi,\omega',s)$ for $\Psi \in \mathscr{S}(Z_{\beta_1\mathbb{A}})$, $\omega'$ a character of $\mathbb{A}^1/k^\times \times \mathrm{GL}(3)^0/\mathrm{GL}(2)_k$, and $s \in \mathbb{C}$. The distribution $\Sigma_{\beta_1}(\Psi,\omega',s)$ is holomorphic at $s = 2$. The stratum $S_{\beta_4}$ corresponds to the prehomogeneous vector space $\mathrm{GL}(2) \times \mathrm{GL}(2)/\widetilde{T}$, $V = k^2 \otimes k^2$, where $\widetilde{T}$ is as before. We define $\Sigma_{\beta_4}(\Psi,\omega',s)$ by (4.4). The distribution $\Sigma_{\beta_4}(\Psi,\omega',s)$ is holomorphic at $s = -2$. The strata $S_{\beta_6}, S_{\beta_8}$ correspond to $\mathrm{GL}(2) \times \mathrm{GL}(1)^2/\widetilde{T}$, $V = \mathrm{Sym}^2 k^2 \oplus k$ or $\mathrm{Sym}^2 k^2 \oplus k^2$ where $\mathrm{GL}(2)$ acts on $k$ trivially, and the action of $\mathrm{GL}(1)^2$ is similar as in §2. A similar argument as in §2 holds for $S_{\beta_6}$ also. We choose $(b_1,b_2) = (3,5),(3,4)$ for $S_{\beta_6}, S_{\beta_8}$ respectively. Let $a_{\beta_6}(\lambda_1) = (1,\lambda_1^5,\lambda_1^{-3}), a_{\beta_8}(\lambda_1) = (1,\lambda_1^4,\lambda_1^{-3})$ for $\lambda_1 \in \mathbb{R}_+$. Then we can choose $\bar{T} = \{a_{\beta_6}(\lambda_1) \mid \lambda_1 \in \mathbb{R}_+\}$ or $\{a_{\beta_8}(\lambda_1) \mid \lambda_1 \in \mathbb{R}_+\}$. We use $d^\times \lambda_1$ on $\bar{T}$. Let $\chi_{\beta_6}, \chi_{\beta_8}$ be principal quasi-characters of $G_\mathbb{A}^1/G_k$ satisfying $\chi_{\beta_6}(a_{\beta_6}(\lambda_1)) = \lambda_1^{-1}$, $\chi_{\beta_8}(a_{\beta_8}(\lambda_1)) = \lambda_1^2$. We define $\Sigma_{\beta_8}(\Psi,\omega',\chi_{\beta_8},s)$ by (4.4). This case was handled by F. Sato in a slightly different formulation [14]. Note that $\omega' = (\omega'_1,\omega'_2,\omega'_3)$ is a character on $\mathrm{GL}(2)^0/\mathrm{GL}(2)_k \times (\mathbb{A}^\times/k^\times)^2$. The case $\mathrm{Sym}^2 k^2 \oplus k$ is of incomplete type, and let $\Sigma_{\beta_6,\mathrm{ad}}(\Psi,\omega',\chi_{\beta_6},s)$ be the adjusted zeta function (see [26] for the details). The distribution $\Sigma_{\beta_6}(\Psi,\omega',\chi_{\beta_6},s)$ is holomorphic at $s = -3$, but $\Sigma_{\beta_8}(\Psi,\omega',\chi_{\beta_8},s)$ has a possible simple pole at $s = -3$. We choose $(x_{1,33},x_{2,13},x_{2,22})$ (resp. $(x_{1,23},x_{2,13},x_{2,22})$) as the coordinate of $Z_{\beta_7}$ (resp. $Z_{\beta_9}$). The stratum $S_{\beta_{10}}$ corresponds to $\mathrm{GL}(2) \times \mathrm{GL}(2)/\widetilde{T}, V = \mathrm{Sym}^2 k^2 \otimes k^2$. We define $\Sigma_{\beta_{10}}(\Psi,\omega',s)$ by (4.4). This case is of incomplete type, and requires an adjusting term. Let $\Sigma_{\beta_{10},\mathrm{ad}}(\Psi,\omega',s)$ be the adjusted zeta function. The distribution $\Sigma_{\beta_{10},\mathrm{ad}}(\Psi,\omega',s)$ is holomorphic at $s = 3$. Let $Z'_{\beta_6,0}$ (resp. $Z'_{\beta_{10},0}$) be the set of points $x$'s such that $x_{i,jk} = 0$ except for $x_{1,33}, x_{2,12}$ (resp. $x_{1,33}, x_{2,22}$). We define

$$R'_{\beta_6,0}\Phi(x_{1,33},x_{2,12}) = \int_\mathbb{A} R_{\beta_6}\Phi(x_{1,33},0,x_{2,12},x_{2,22})dx_{2,22},$$

$$R'_{\beta_{10},0}\Phi(x_{1,33},x_{2,22}) = \int_{\mathbb{A}^2} R_{\beta_{10}}\Phi(0,0,x_{1,33},x_{2,22},x_{2,23},x_{2,33})dx_{2,23}dx_{2,33}.$$

We need some definitions in order to state the result in this case.

**Definition (8.6)**

(1) $\delta_{\beta_1}(\omega) = \delta(\omega_2)$, $\omega_{\beta_1} = (1, \omega_1)$,
(2) $\delta_{\beta_4}(\omega) = \delta(\omega_1)$, $\omega_{\beta_4} = (1, \omega_2)$,
(3) $\delta_{\beta_6}(\omega) = \delta(\omega_1^2)$, $\omega_{\beta_6} = (\omega_1, \omega_2, \omega_2)$,
(4) $\delta_{\beta_7}(\omega) = \delta(\omega_1^3)\delta(\omega_2^2)$, $\omega_{\beta_7} = (\omega_2, \omega_1, \omega_1^{-1}\omega_2)$,
(5) $\delta_{\beta_8}(\omega) = \delta(\omega_1)$, $\omega_{\beta_8} = (\omega_1, \omega_2, \omega_2)$,
(6) $\delta_{\beta_9}(\omega) = \delta(\omega_1^3)\delta(\omega_2)$, $\omega_{\beta_9} = (1, \omega_1^{-1}, \omega_1)$,
(7) $\delta_{\beta_{10}}(\omega) = \delta(\omega_1)$, $\delta_{\beta_{10},st}(\omega) = \delta(\omega_1)\delta(\omega_2^2)$, $\omega_{\beta_{10}} = (1, \omega_2)$, $\omega_{\beta_{10},st} = (\omega_2, \omega_2)$.

We define

$$F_{-1}(\Phi, \omega, 0) = \mathfrak{V}_2 \mathfrak{V}_3 \delta_\#(\omega) \Phi(0) + 2\delta_{\beta_4}(\omega) \Sigma_{\beta_4}(R_{\beta_4}\Phi, \omega_{\beta_4}, -2)$$
$$+ 8\delta_{\beta_6}(\omega) \Sigma_{\beta_6,\mathrm{ad}}(R_{\beta_6}\Phi, \omega_{\beta_6}, \chi_{\beta_6}, -3)$$
$$+ 14\delta_{\beta_8}(\omega) \Sigma_{\beta_8,(0)}(R_{\beta_8}\Phi, \omega_{\beta_8}, \chi_{\beta_8}, -3))$$
$$+ \frac{\delta_\#(\omega)}{2} \left( 5\Sigma_{2,(1,0)}(R'_{\beta_6,0}\Phi, -1, -3) + 3\Sigma_{2,(0,1)}(R'_{\beta_6,0}\Phi, -1, -3) \right),$$

$$F_{-2}(\Phi, \omega, 0) = -3\delta_\#(\omega) \Sigma_2(R'_{\beta_6,0}\Phi, -1, -3),$$

$$F_{-1}(\Phi, \omega, 2) = \delta_{\beta_1}(\omega) \Sigma_{\beta_1}(R_{\beta_1}\Phi, \omega_{\beta_1}, 2) + \frac{\delta_{\beta_9}(\omega)}{3} \Sigma_3(R_{\beta_9}\Phi, \omega_{\beta_9}, -1, -\frac{2}{3}, \frac{2}{3}),$$

$$F_{-1}(\Phi, \omega, 3) = \frac{\delta_{\beta_7}(\omega)}{6} \Sigma_{3,(0,0,0)}(R_{\beta_7}\Phi, \omega_{\beta_7}, \frac{1}{2}, 0, \frac{1}{2})$$
$$+ \delta_{\beta_{10}}(\omega) \Sigma_{\beta_{10},\mathrm{ad}}(R_{\beta_{10}}\Phi, \omega_{\beta_{10}}, 3)$$
$$- \frac{\delta_{\beta_{10},st}(\omega)}{8} \Sigma_{2,(1,0)}(R'_{\beta_{10},0}\Phi, \omega_{\beta_{10},st}, \frac{1}{2}, \frac{1}{2})$$
$$+ \frac{\delta_{\beta_{10},st}(\omega)}{24} \Sigma_{2,(0,1)}(R'_{\beta_{10},0}\Phi, \omega_{\beta_{10},st}, \frac{1}{2}, \frac{1}{2}),$$

$$F_{-2}(\Phi, \omega, 3) = -\frac{\delta_{\beta_{10},st}(\omega)}{4} \Sigma_2(R'_{\beta_{10},0}\Phi, \omega_{\beta_{10},st}, \frac{1}{2}, \frac{1}{2}),$$

$$F(\Phi, \omega) = \sum_{i=0,3, j=1,2} \frac{F_{-j}(\Phi, \omega, i)}{(s-i)^j} + \frac{F_{-1}(\Phi, \omega, 2)}{s-2}.$$

We proved the following theorem in [26].

**Theorem (8.7)** *Suppose that $M_{V,\omega}\Phi = \Phi$. Then*

$$Z(\Phi, \omega, s) = Z_+(\Phi, \omega, s) + Z_+(\widehat{\Phi}, \omega^{-1}, 12-s) - F(\widehat{\Phi}, \omega^{-1}, 12-s) - F(\Phi, \omega, s).$$

Most of the coefficients of the Laurent expansions at the poles of $Z(\Phi, \omega, s)$ are generalized special values of the adjusted zeta functions defined for $Z_\beta$'s. However, we have some extra terms like $\Sigma_{2,(1,0)}(R'_{\beta_6,0}\Phi, -1, -3)$. The author expects similar phenomena for other prehomogeneous vector spaces. Also the author hopes to handle more cases in the near future.

## References

[1] Atiyah, M. F. Convexity and commuting hamiltonians. *Bull. London Math. Soc.*, 14:1–15, 1982.

[2] Blanchard, A. Sur les variétés abalytiques complexes. *Ann. Sci. Ecole Norm. Sup.*, 73:178–201, 1956.

[3] Datskovsky, B., and D. J. Wright. Density of discriminants of cubic extensions. *J. Reine Angew. Math.*, 386:116–138, 1988.

[4] Davenport, H. On the class-number of binary cubic forms I and II. *London Math. Soc.*, 26:183–198, 1951. Corrigendum: ibid., 27:512, 1952.

[5] Davenport, H., and H. Heilbronn. On the density of discriminants of cubic fields I. *Bull. London Math. Soc.*, 1:345–348, 1961.

[6] Davenport, H. and H. Heilbronn. On the density of discriminants of cubic fields. II. *Proc. Royal Soc.*, A322,:405–420, 1971.

[7] Guillemin, V., and S. Sternberg. Convexity properties of the moment mapping. *Invent. Math.*, 67:491–513, 1982.

[8] Kempf, G. Instability in invariant theory. *Ann. of Math.*, 108:299–316, 1978.

[9] Kempf, G., and L. Ness. The length of vectors in representation spaces. In *Algebraic Geometry, Proceedings, Copenhagen*, volume 732, pages 233–242. Springer Lecture Notes in Mathematics, Berlin, Heidelberg, New York, 1978.

[10] Kirwan, F. C. *Cohomology of quotients in symplectic and algebraic geometry*. Mathematical Notes, Princeton University Press, 1984.

[11] Langlands, R. The volume of the fundamental domain for some arithmetical subgroups of Chevalley groups. In *Proceedings of symposium in pure mathematics*, volume 9. 1966.

[12] Mumford, D., and J. Fogarty. *Geometric invariant theory, 2nd edition*. Springer, Berlin, Heidelberg, New York, 1982.

[13] Ness, L. A stratification of the null cone via the moment map. *Amer. J. Math.*, 106:1281–1325, 1984.

[14] Sato, F. On zeta functions of ternary zero forms. *J. Fac. Sci. Univ. Tokyo, Sect IA*, 28:585–604, 1981.

[15] Sato, F. Zeta functions in several variables associated with prehomogeneous vector spaces I: Functional equations. *Tôhoku Math. J.*, (2) 34:no. 3 437–483, 1982.

[16] Sato, F. Zeta functions in several variables associated with prehomogeneous vector spaces II: A convergence criterion. *Tôhoku Math. J.*, (2) 35 no. 1:77–99, 1983.

[17] Sato, M., and T. Shintani. On zeta functions associated with prehomogeneous vector spaces. *Ann. of Math.*, 100:131–170, 1974.

[18] Shintani, T. On Dirichlet series whose coefficients are class-numbers of integral binary cubic forms. *J. Math. Soc. Japan*, 24:132–188, 1972.

[19] Shintani, T. On zeta-functions associated with vector spaces of quadratic forms. *J. Fac. Sci. Univ. Tokyo, Sect IA*, 22:25–66, 1975.

[20] Siegel, C. L. The average measure of quadratic forms with given discriminant and signature. *Ann. of Math.*, 45:667–685, 1944.

[21] Weil, A. *Basic number theory*. Springer, Berlin, Heidelberg, New York, 1974.

[22] Wright, D. J. The adelic zeta function associated to the space of binary cubic forms part I: Global theory. *Math. Ann.*, 270:503–534, 1985.

[23] Wright, D. J., and A. Yukie. Prehomogeneous vector spaces and field extensions. *Invent. Math.*, 110:283–314, 1992.

[24] Ying, K. *On the convergence of the adelic zeta function associated to irreducible regular prehomogeneous vector spaces*. Ph. D. Thesis, Johns Hopkins University, 1993.

[25] Yukie, A. On the Shintani zeta function for the space of binary quadratic forms. *Math. Ann.*, 292:355–374, 1992.

[26] Yukie, A. *Shintani zeta functions*, volume 183. London Mathematical Society Lecture Note series, Cambridge, 1993.

Akihiko Yukie
Oklahoma State University
Mathematics Department
401 Math Science
Stillwater OK 74078-1058 USA
yukie@math.okstate.edu